\documentclass[10pt,letterpaper,conference]{ieeeconf}

\usepackage{latexsym, amssymb, amsmath}
\usepackage[dvips]{graphicx} % added for inserting the pictures

\usepackage{float}      %these two packages are for fixing figures in place
\restylefloat{figure}
\usepackage{graphicx,dsfont}
\usepackage{mathrsfs}  % package for math fonts: mathscripts \mathscr
\usepackage{color} % standard color package
\usepackage[latin1]{inputenc}
\usepackage{enumerate}
\usepackage{flushend} % align at the last lines of every page
\usepackage{stfloats}
\usepackage{cite} % package to reduce the length of long citations by using a hyphen

\usepackage{balance}

\usepackage{multicol}
\usepackage[german,english]{babel}
\usepackage[T1]{fontenc}
\usepackage[latin1]{inputenc}

\allowdisplaybreaks[4] % added by Hong for allowing page breaks among a large formula

\newtheorem{prop}{Proposition} % [section]
\newtheorem{proposition}[prop]{Proposition} % [section]
\newtheorem{lem}{Lemma}%[section]
\newtheorem{lemma}[lem]{Lemma} % [section]
\newtheorem{thm}{Theorem} % [section]
\newtheorem{theorem}[thm]{Theorem} % [section]
\newtheorem{cor}{Corollary} % [section]
\newtheorem{corollary}[cor]{Corollary} % [section]
\newtheorem{defn}{Definition} % [section]
\newtheorem{definition}[defn]{Definition} % [section]
\newtheorem{exmp}{Example} % [section]
\newtheorem{example}[exmp]{Example} % [section]
\newtheorem{exam*}{Example}

\def\custombibliography#1{
 \normalsize
\section*{\centering References}
 \list
 {[\arabic{enumi}]}{\settowidth\labelwidth{[#1]}\leftmargin\labelwidth
 \setlength{\itemsep}{.1em}
 \advance\leftmargin\labelsep
 \usecounter{enumi}}
 \def\newblock{\hskip .11em plus .33em minus -.07em}
 \sloppy
 \sfcode`\.=1000\relax}

\def\L2{{\cal L}_2}
\def\bull{\rule{0.08in}{0.08in}} % square filled bullet
\def\openbull{\framebox[0.08in][c]{$\;$}} % square unfilled bullet
\def\re{{\mathbb R}} % real numbers (AMS symbol)

\def\shuffle{{\scriptscriptstyle \;\sqcup \hspace*{-0.05cm}\sqcup\;}}

\def\allseriesell{\mbox{$\re^{\ell} \langle\langle X \rangle\rangle $}}

\def\eqref#1{(\ref{#1})} % parentheses around referenced equation numbers
\def\abs#1{\left\vert #1 \right\vert}
\def\norm#1{\left\Vert#1\right\Vert}

\def\supp{{\rm supp}}

\newcommand{\comment}[1]{} % Allows one to comment out a block of text

\def\begce{\begin{center}}
\def\endce{\end{center}}
\def\begar{\begin{array}}
\def\endar{\end{array}}
\def\begeq{\begin{equation}}
\def\endeq{\end{equation}}
\def\begdi{\begin{displaymath}}
\def\enddi{\end{displaymath}}
\def\begdis{\begin{eqnarray*}}
\def\enddis{\end{eqnarray*}}
\def\begeqa{\begin{eqnarray}}
\def\endeqa{\end{eqnarray}}
\def\begdes{\begin{description}}
\def\enddes{\end{description}}
\def\begit{\begin{itemize}}
\def\endit{\end{itemize}}
\def\begen{\begin{enumerate}}
\def\enden{\end{enumerate}}
\def\beglar{\left[\begin{array}}
\def\endrar{\end{array}\right]}
\def\begle{\begin{lemma}}
\def\endle{\end{lemma}}
\def\begde{\begin{definition}}
\def\endde{\end{definition}}
\def\begth{\begin{theorem}}
\def\endth{\end{theorem}}
\def\begco{\begin{corollary}}
\def\endco{\end{corollary}}
\def\begprop{\begin{proposition}}
\def\endprop{\end{proposition}}
\def\begex{\begin{example}}
\def\endex{\hfill\openbull \end{example} \vspace*{0.1in}}
\def\begexer{\begin{exercise}}
\def\endexer{\end{exercise}}

\def\begres{\noindent{\bf Remarks}:\begin{enumerate}}
\def\endres{\end{enumerate} \par}
\def\begpr{\noindent{\em Proof:}$\;\;$}
\def\endpr{\hfill\bull \vspace*{0.1in}}

\def\begtab{\begin{tabular}}
\def\endtab{\end{tabular}}
\def\rref#1{(\ref{#1})}
\newcommand\cdcout[1]{} % reduce it to 6 pages
\newcommand{\rv}[1]{\boldsymbol{#1}} % Use italic boldface to indicate the Random Variables
\newcommand{\RomanNumber}[1]{\uppercase\expandafter{\romannumeral #1}}
\newcommand{\romannumber}[1]{\lowercase\expandafter{\romannumeral #1}}

\DeclareMathAlphabet{\mathpzc}{OT1}{pzc}{m}{it}

\def\1{\rv 1} %Indicator
\def\allseriesPXelln{\mbox{$\re^{\ell\times n}\langle\langle 
\mathfrak{P}X\rangle\rangle$} }

\def\allseriesTXelln{\mbox{$\re^{\ell\times n}\langle\langle 
\mathfrak{T}X\rangle\rangle$} }
\def\allpolyTXelln{\mbox{$\re^{\ell\times n}\langle 
\mathfrak{T}X\rangle$} }
\def\allseriesTXnn{\mbox{$\re^{n\times n}\langle\langle 
\mathfrak{T}X\rangle\rangle$} }

\def\allseriesTDXelln{\mbox{$\re^{\ell\times n}\langle\langle 
\mathfrak{TD}X\rangle\rangle$} }
\def\allpolyTDX{\mbox{$\re\langle 
\mathfrak{TD}X\rangle$} }
\def\allpolyTDXelln{\mbox{$\re^{\ell\times n}\langle 
\mathfrak{TD}X\rangle$} }
\def\allpolyTXnn{\mbox{$\re^{n\times n}\langle 
\mathfrak{T}X\rangle$} }
\def\allpolyTX{\mbox{$\re\langle \mathfrak{T}X\rangle$} }

\def\allseriesA{\mbox{$\re\langle\langle A\rangle\rangle$}}
\def\allseriesA'{\mbox{$\re\langle\langle A'\rangle\rangle$}}

\def\allpolyell{\mbox{$\re^{\ell}\langle X\rangle$}}

\def\L1spaceprodu{{ L}_1(\Omega\times [0,T],{\mathcal P},P\otimes \lambda)}

\def\Hspace0{{\mathcal H}^2_0}

\def\shuffleNC{{ \prec \hspace*{-0.07cm}\succ}}
\newcommand{\arb}[1]{\begin{matrix}\includegraphics[height=7mm]{a#1.eps}
\end{matrix}}
\newcommand{\arblarge}[1]{\begin{matrix}\includegraphics[height=9mm]{a#1.eps}
\end{matrix}}
\newcommand{\arbsize}[1]{\begin{matrix}\includegraphics[height=14mm]{a#1.eps}
\end{matrix}}                      
\newcommand{\arbsmall}[1]{\begin{matrix}\includegraphics[height=3mm]{a#1.eps}
\end{matrix}}                                            
                      
\newcount\colveccount
\newcommand*\colvec[1]{
        \global\colveccount#1
        \begin{pmatrix}
        \colvecnext
}
\def\colvecnext#1{
        #1
        \global\advance\colveccount-1
        \ifnum\colveccount>0
                \\
                \expandafter\colvecnext
        \else
                \end{pmatrix}
        \fi
}

\renewcommand{\baselinestretch}{0.99}

\allowdisplaybreaks[4]

\IEEEoverridecommandlockouts                    % This command is only
\title{Dendriform-Tree Setting for Fully Non-commutative Fliess 
Operators$^\ast$\thanks{$^\ast$The second author was supported by a 
grant from the Severo Ochoa Excellence Program at the Instituto de Ciencias 
Matem\'{a}ticas in Madrid, Spain. The third author was supported by
a Ram\'on y Cajal research grant from the Spanish government.}}

\author{Luis A. Duffaut Espinosa$^\dag$\thanks{$^\dag$Department of Electrical 
and Computer Engineering, George Mason University, Fairfax, VA 22030, 
USA, {\tt\small lduffaut@gmu.edu}} $\;\;\;$ W. Steven Gray$^{\ddag 
\,\S}$\thanks{$^\ddag$Instituto de Ciencias Matem\'{a}ticas, Consejo 
Superior de Investigaciones Cient\'ificas, 28049 Madrid, Spain {\tt\small 
sgray@odu.edu, kurusch@icmat.es}} $\;\;\;$ Kurusch 
Ebrahimi-Fard$^\ddag$\thanks{$^\S$On leave from Old Dominion University, 
Norfolk, VA 23529, USA}}

\begin{document}

\maketitle

\begin{abstract}
This paper provides a dendriform-tree setting for Fliess operators with 
matrix-valued inputs. This class of analytic nonlinear input-output systems is 
convenient, for example, in quantum control. In particular, a 
description of such Fliess operators is provided using planar binary trees. 
Sufficient conditions for convergence of the defining series are also 
given.
\end{abstract}

\section{Introduction}  \label{sec:1}

Fliess operators provide a general framework under which analytic nonlinear 
input-output systems can be studied 
\cite{Fliess_81,Gray-Wang_SCL02,Gray-Li_05}. Let $X=\{x_0,x_1,\ldots,x_m\}$ be 
an alphabet and $X^{\ast}$ the free monoid comprised of all words over $X$ 
(including the empty word $\emptyset$) under the catenation product. A formal 
power series $c$ in $X$ is any mapping of the form $X^{\ast}\rightarrow 
\re^{\ell}\!:\eta\mapsto (c,\eta)$. The set of all such mappings will be 
denoted by $\allseriesell$. The support of an arbitrary series $c$ is 
$\supp(c)=\{\eta\in X^\ast, (c,\eta) \neq 0\}$.  A series having finite support 
is called a \emph{polynomial}, and the set of all polynomials is 
represented by $\allpolyell$. For a measurable function $u: [a, b] 
\rightarrow\re^m$ define $\norm{u}_{L_p}=\max\{\norm{u_i}_{L_p}: \ 1\le i\le 
m\}$, where $\norm{u_i}_{L_p}$ is the usual $L_p$-norm for a measurable 
real-valued component function $u_i$. Define iteratively for each $\eta\in 
X^{\ast}$ the mapping $E_\eta: L_1^m[t_0, t_0+T]\rightarrow C[t_0, t_0+T]$ by 
$E_\emptyset[u] = 1$, and
\begeq \label{eq:comm_E_map}
E_{x_i{\eta'}}[u](t,t_0) = 
\int_{t_0}^tu_{i}(\tau)E_{{\eta'}}[u](\tau,t_0)\,d\tau,
\endeq
where $x_i\in X$, ${\eta'}\in X^{\ast}$ and $u_0= 1$. The input-output 
operator corresponding to $c$ is then
\begdi
F_c[u](t) :=  \sum_{\eta\in X^{\ast}}
(c,\eta)\,E_\eta[u](t),
\enddi
which is called a \emph{Fliess operator}. If the generating series $c$ is {\em 
locally convergent}, i.e., there exists constants $K,M>0$ such that 
%\begdi 
$\abs{(c,\eta)} \,\le \,K M^{|\eta|}|\eta|!$ 
%\enddi
for all $\eta\in X^{\ast}$, where $\abs{\eta}$ denotes the number of letters
in $\eta$, then $F_c[u]$ converges absolutely and uniformly on $[t_0,t_0+T]$ 
if $T$ and $\norm{u}_{L_p}$ are sufficiently small. In general, the 
input-output map $F_c:u\rightarrow y$ need not have a state space realization, 
however, many familiar and relevant examples are obtained from the state space 
setting.

A tacit assumption in the standard theory of Fliess operators is that the 
inputs are mutually commutative, i.e., the functions associated with each 
letter of $X$ commute pointwise in time. The proposition 
here is that this assumption results in a great deal of simplification and 
hides certain underlying algebraic structures that are important in 
applications like control on Lie groups \cite{Brockett_73} and quantum control 
\cite{Alessandro_2007}.

As a motivating example, consider a bilinear system 
\begin{equation} \label{eq:motiv_example}
\dot{z}(t)= A z(t)+B(t)z(t)u(t),
\end{equation}
where $B$ is a smooth function on $[0,T]$. One can view $u$ as the user 
controlled input and $B$ as a disturbance input. Let $z_i$ be the solution of 
\rref{eq:motiv_example} when $z(0)=e_i=[0,\hdots,0,1,0,\hdots,0]^T$ with the $1$ 
in the $i$-th position and define $Z(t)= [x_1(t), \cdots, x_n(t)]$, where $n$ is 
the dimension of the system. Then
\begin{align} \label{eq:operator_evolution}
 \dot{Z}(t)  & = (A +B(t)u(t))Z(t) =: U(t)Z(t), 
\end{align}
where in general $U(t_1)U(t_2) \neq U(t_2)U(t_1)$. This is, for example, the 
setting of a regulator problem in which the input-output map from disturbance to 
some output $y(t)=CZ(t)$ needs to be determined when $u(t)=u_0\in \re$. 
Equation \rref{eq:operator_evolution} is also the usual starting point for 
control theory on Lie groups. 
% The product $U(t)Z(t)$ is viewed as a 
% left-invariant vector field. 
Systems such as in \rref{eq:operator_evolution} 
are ubiquitous in quantum mechanics. Take for instance the case of a spin 
particle in a magnetic field $B_m$ whose direction changes in time. The function 
$U$ is proportional to the scalar product $S \cdot B_m$, where $S$ represents 
the spin vector. Now suppose the magnetic field at $t = t_1$ is parallel to the 
$x$-axis, and at $t=t_2$ to the $y$-axis, then $U(t_1) \varpropto \abs{B_m} 
S_x$, $U(t_2) \varpropto \abs{B_m} S_y$, and $[U(t_1), U(t_2)] \varpropto B_m^2 
[S_x , S_y ] \varpropto B_m^2 S_z \neq 0$. Moreover, systems of the form 
$\dot{Z}= U(t)F(Z(t))$ can be considered where the coordinate change 
$\bar{Z}=F(Z)$ is valid on a neighborhood of $Z(0)=I$. In which case,
\begdi
\dot{\bar{Z}}(t) = 
\left(\frac{dF^{-1}(\bar{Z})}{d\bar{Z}}\left.\rule{0in}{0.2in}\right|_{\bar{Z}
=\bar{Z} (0)} \right)^ {-1} U(t)\bar{Z}(t) =: W(t)\bar{Z}(t),
\enddi
is in the same class as \rref{eq:operator_evolution}. 

A series representation of the solution of \rref{eq:operator_evolution} can 
be obtained by successive iterations. That is, 
\begin{align} \label{eq:dyson_series}
Z(t)= & {\,} I+ \sum_{n=1}^\infty \int_{0}^{t}\!\! U(t_1) dt_1 
\cdots \int_{0}^{t_{n-1}} \!\!\!\!\!\!\!\!\!\!U(t_n)dt_n.
\end{align}
This series has an artificial exponential representation in terms of the 
\emph{time ordered operator} 
\begin{align*}
 \mathcal{T} \!\!\left( U(t_1) \cdots 
U(t_n)\right) := & \sum_{\sigma\in S_n}\Theta_n^\sigma 
U(t_{\sigma(1)}) \cdots U(t_{\sigma(n)}),
\end{align*}
where $\Theta_n^\sigma = 
\prod_{i=1}^{n-1}\Theta(t_{\sigma(i)}-t_{\sigma(i+1)})$, $\Theta$ is the 
Heaviside step function, $\sigma$ is a permutation, and $S_n$ is the group of 
all permutations of order $n$ \cite{Bauer}. Because of the symmetry of the 
simplex consisting of all ordered $n$-tuples 
$(t_1,t_2,\cdots ,t_n)$ in the integration limits, this operator 
satisfies:
\begin{align*}
\lefteqn{\int_{t_0}^{t} dt_1 \int_{t_0}^{t_1} dt_2 \cdots \int_{t_0}^{t_{n-1}} 
dt_n\, U(t_1)U(t_2) \cdots U(t_n)} & \\
 & = \frac{1}{n!}\int_{t_0}^{t}\!\! dt_1 \int_{t_0}^{t}\!\!\!\! dt_2 
\cdots \int_{t_0}^{t}\!\!\! dt_n \,\mathcal{T}\!\left(U(t_1)U(t_2) \cdots 
U(t_n)\right).
\end{align*}
The solution is thus written as the time ordered exponential 
\begin{align} \label{eq:time_ordered_exp}
%\nonumber \lefteqn{X(t)}\\ 
\nonumber    Z(t)  & = I+ \sum_{n=1}^\infty \frac{1}{n!} \int_{t_0}^{t}
  \!\!\cdots\! \int_{t_0}^{t} \mathcal{T}\left(U(t_1) \cdots U(t_n)\right)dt_1 
  \cdots dt_n\\
     & =: \mathcal{T}\exp \left( \int_{t_0}^{t} U(s)\, ds \right).
\end{align}
Expression \rref{eq:time_ordered_exp} disregards the algebra 
provided by the non-commutative iterated integrals in \rref{eq:dyson_series}. 
However, it is known that by systematically keeping track of the 
non-commutative orderings of the iterated integrals, a true exponential 
(\emph{Magnus expansion}) can be derived. That is, $X(t)  = \exp 
\left(\Omega(U(t)) \right)$, where $\Omega$ is obtained via a recursion 
\cite{Kurusch-Manchon_2009,Kurusch-Manchon_2012,Magnus_54 }. In the case of 
commutative inputs, the algebra provided by the iterated integrals is the 
\emph{shuffle algebra}, which is based on the integration by parts formula 
\cite{Reutenauer_93,Ree_57}. The noncommutative version of this formula is 
\begin{align*}
\int_0^t u_i(s) \,ds \int_0^t u_j(s) \,ds = & {\,} \int_0^t u_i(s) 
\left(\int_0^s u_j(r) \,dr \right)ds \\
& + \int_0^t \left(\int_0^s u_i(r)  \,dr \right) u_j(s) \,ds. 
\end{align*}
Note that the second summand on the right-hand side above cannot be generated 
recursively as in \rref{eq:comm_E_map}. Moreover, products of iterated 
integrals are fundamental when the system's state is filtered by an analytic 
output function \cite{Fliess_81,Wang_90}, in the computation of bounds 
for iterated integrals \cite{Duffaut-Gray-Gonzalez_2012} and the 
characterization of system interconnections such as the product, cascade and 
feedback connections \cite{Gray-Li_05}. The first goal of this paper is to 
provide a fully non-commutative extension of the theory of Fliess operators in 
the context of \emph{dendriform/tree} algebras. They will be referred to as 
\emph{dendriform Fliess operators}. The second goal is to give sufficient 
conditions under which dendriform Fliess operators with non-commutative inputs 
converge.

The paper is organized as follows. Section \ref{sec:2} provides a 
tutorial treatment of dendriform algebras. In Section \ref{sec:3}, planar 
binary trees are presented as a tool to keep track of the non-commutativity of 
iterated integrals. Also, the non-commutative version of the shuffle product is 
given. These results are then applied in Section \ref{sec:4} to define 
dendriform Fliess operators. Then sufficient conditions for the convergence 
of dendriform Fliess operators are provided. Finally, the conclusions are given 
in Section \ref{sec:conclusions}.

\section{Dendriform Algebras} \label{sec:2}

The goal of this section is to introduce parenthesis words and their 
relationship to dendriform algebras. The concepts here can be found in 
\cite{Loday-Ronco_98,Kurusch-Guo_04} and references therein.

Let $X$ be a finite alphabet and $\mathfrak{P}X=X\cup 
\{\lfloor,\rfloor\}$. The free semigroup under catenation generated by 
$\mathfrak{P}X$ is denoted $\mathfrak{P}X'$. For $\eta=q_{1}q_{2}\cdots 
q_{n}\in \mathfrak{P}X'$, let $s(\eta)_i$ denote the number of $\lfloor$'s in 
$q_{1}\cdots q_{i}$ minus the number of $\rfloor$'s in $q_{1}\cdots q_{i}$. 
\begde \label{de:parenthesis_words} A word $\eta=q_{1}q_{2}\cdots q_{n} \in 
\mathfrak{P}X'$ is called a \emph{parenthesis word} if its parenthesization is 
balanced, i.e., it satisfies:
\begin{itemize} 
 \item[$i.$] $s(\eta)_i \ge 0$ for $i=1,\hdots,n-1$ and $s(\eta)_n=0.$
 \item[$ii.$] $q_{i}q_{i+1}\neq x_{i_1}x_{i_2}$ for $x_{i_1},x_{i_2}\in X$ and  
$i=1,\hdots,n-1$. 
 \item[$iii.$] $q_{i}q_{i+1}\neq \lfloor\rfloor, \rfloor\lfloor$ for 
$i=1,\hdots,n-1$.
 \item[$iv.$] $q_1=\lfloor$ and $q_n=\rfloor$ cannot occur at the same time.
 \item[$v.$] There are no sub-words in $\eta$ of the form $\xi \lfloor \nu 
\rfloor \kappa$ or $\lfloor \lfloor \xi \rfloor  \rfloor$ for 
$\xi,\nu,\kappa\in \mathfrak{P}X'$.
\end{itemize}
\endde
Parenthesis words are such that $x_i\lfloor x_j \rfloor \neq \lfloor x_i\rfloor 
x_j$ for $x_i,x_j\in X$. This set of parenthesis words constitutes a \emph{free 
Magma} under balanced parenthesization 
\cite{Melancon_92,Kurusch-Guo_04,Kurusch-Guo_08}. The set of 
parenthesis words including the empty word $\emptyset$ is denoted by 
$\mathfrak{P}X^\ast$. In Section \ref{sec:3}, the operation in this magma is 
better understood in terms of the grafting operation on trees. A formal 
power series in $\mathfrak{P}X$ is any mapping of the 
form $\mathfrak{P}X^\ast\rightarrow \re^{\ell\times n}\!:\eta\mapsto (c,\eta)$. 
The set of all such mappings will be denoted by $\allseriesPXelln$, which forms 
an $\re$-vector space.

An alternative to parenthesization of words is to encode the order in which 
balanced parentheses appear by using two different products, say $\prec$ and 
$\succ$. For example, 
\begeq \label{eq:parenthesis_to_dendriform1}
x_i\lfloor x_j\rfloor\equiv x_i\prec x_j \mbox{ and } \lfloor x_i\rfloor x_j 
\equiv x_i \succ x_j.
\endeq
Using these products the induced algebraic structure on $\mathfrak{P}X^\ast$ is 
described next. 

\begde
A \emph{dendriform algebra} is an $\re$-vector space, $(D,+,\cdot)$, endowed 
with products $\prec$ and $\succ$ such that for $a,b,c\in D$ the following 
axioms are satisfied:
\begin{subequations} \label{subeq:dendriform_axioms}
\begin{align}
\label{eq:Dend_axiom1} (a \prec b) \prec c  & =  a \prec (b \prec c 
+ b \succ c) , \\
\label{eq:Dend_axiom2} (a \succ b) \prec c  & = a \succ ( b \prec c 
) , \\
\label{eq:Dend_axiom3} a \succ ( b \succ c ) & =  ( a \prec b + a 
\succ b ) \succ c.
\end{align}
\end{subequations}
\endde
If $D=X$, then $(X,\prec,\succ)$ forms a dendriform algebra. Similar to 
\rref{eq:parenthesis_to_dendriform1}, for every $\eta\in \mathfrak{P}X^\ast$ 
there is a corresponding dendriform product in $(X,\prec,\succ)$. This is made 
explicit by the injection $\delta: \mathfrak{P}X^\ast \rightarrow 
(X,\prec,\succ)$, which is defined recursively by
\begin{align*}
\delta(\eta) = \left\{\begin{array}{ll}
                x_i \prec \delta(\eta'), & \mbox{if} \quad \eta = x_i \lfloor 
		\eta' \rfloor, \\
		\delta(\eta') \succ x_i , & \mbox{if} \quad \eta = \lfloor 
		\eta' \rfloor x_i, \\
		\delta(\eta') \succ x_i \prec \delta(\eta'') , & \mbox{if} 
		\quad \eta = \lfloor \eta' \rfloor x_i \lfloor \eta'' \rfloor,\\
               \end{array}\right.
\end{align*}
where $x_i \in X$, $\eta',\eta''\in \mathfrak{P}X^\ast$, 
$\delta(\emptyset)=\emptyset$, and $\delta(x_j)=x_j$ for all $x_j\in X$. For 
example,
\begdi
\delta(x_i\lfloor \lfloor x_j\rfloor x_k\rfloor) = x_i \prec (\delta(  
\lfloor x_j\rfloor x_k) ) = x_i \prec ( x_j \succ x_k ).
\enddi
Define $\mathfrak{T}X^\ast =\delta(\mathfrak{P}X^\ast)$, the image of 
$\mathfrak{P}X^\ast$ under $\delta$. Any element of $\mathfrak{T}X^\ast$ is 
called a \emph{dendriform word}. 

The set of formal power series on dendriform words is denoted 
by $\allseriesTXelln$, and it is also an $\re$-vector space. An element of 
$\allseriesTXelln$ can be viewed as a mapping $c:\mathfrak{T}X^\ast\rightarrow 
\re^{\ell\times n}\!:\eta\mapsto (c,\eta)$. The set of all series in 
$\allseriesTXelln$ having finite support is denoted  by $\allpolyTXelln$. 
In addition, for any dendriform word there is only one corresponding word in 
$X^\ast$ given by the projection $\varphi: \mathfrak{T}X^\ast \rightarrow 
X^\ast$. For example, $\varphi(x_i \prec( x_j \prec x_k )) = x_ix_jx_k \in 
X^\ast$.

Next define the product $\shuffleNC : \mathfrak{T}X^\ast \times 
\mathfrak{T}X^\ast \rightarrow \allpolyTXelln : (\eta,\xi) \mapsto \eta \prec 
\xi + \eta \succ \xi$. This product is the non-commutative counterpart of the 
\emph{shuffle product} \cite{Fliess_81}, and it is extended 
bilinearly on $\allseriesTXelln \times \allseriesTXelln$.

\begle \cite{Loday-Ronco_98,Melancon_92,Kurusch-Manchon_2012} 
\label{le:associativity_shuffleNC} 
$(\allseriesTXelln,\shuffleNC)$ is an associative $\re$-algebra.
\endle

An important characteristic of the commutative shuffle product is that it can 
be defined recursively, which is convenient for computer implementations. For 
the non-commutative shuffle product such a recursive definition is only 
available when the words to be \emph{shuffled} have single letters. In this 
regard, the notion of planar binary trees plays a key role as described next.

\section{Trees, Dendriform Words and Iterated Integrals} \label{sec:3}

The objective of this section is to describe the one-to-one correspondence 
between planar binary trees and dendriform words. Then their relationship to 
non-commutative iterated integrals is described. The majority of concepts 
presented in subsection \ref{subsec:III-A} can be found in 
\cite{Loday-Ronco_98,Melancon_92,Kurusch-Manchon_2012} and references therein.

\subsection{Trees and dendriform words} \label{subsec:III-A}

A tree is a non-cyclic connected graph $(V,\Gamma)$, where $V$ denotes the 
vertices of the graph and $\Gamma$ the edges. A \emph{leaf} is defined as a
vertex that is the endpoint of only one edge. The $n$ leaves of a tree are 
labeled from left to right as $1,2,\hdots,n$. A \emph{planar rooted} tree is 
a tree embedded in the plane in which one vertex (with no incoming edges) is 
labeled as the \emph{root}. The \emph{interior} vertices of a rooted planar 
tree is the set $V$ minus the root and the leaves. A planar \emph{$n$-ary} tree 
is a planar rooted tree where every interior vertex has one root and $n$ 
leaves. \emph{Order} is defined by the number of interior vertices. This paper 
is concerned with \emph{planar binary trees}, so every interior vertex has one 
root and two leaves. The set of all planar binary trees is denoted by 
$\mathfrak{T}$, and $\mathfrak{T}_n$ denotes the set of planar 
binary trees of order $n$. The planar binary trees up to order three are:  
\begin{align*}
\mathfrak{T}_0 &= \left\{ \arb{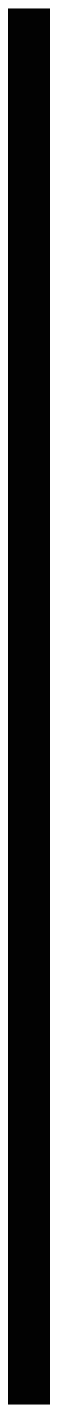} \right\},
\mathfrak{T}_1 = \left\{ \arb{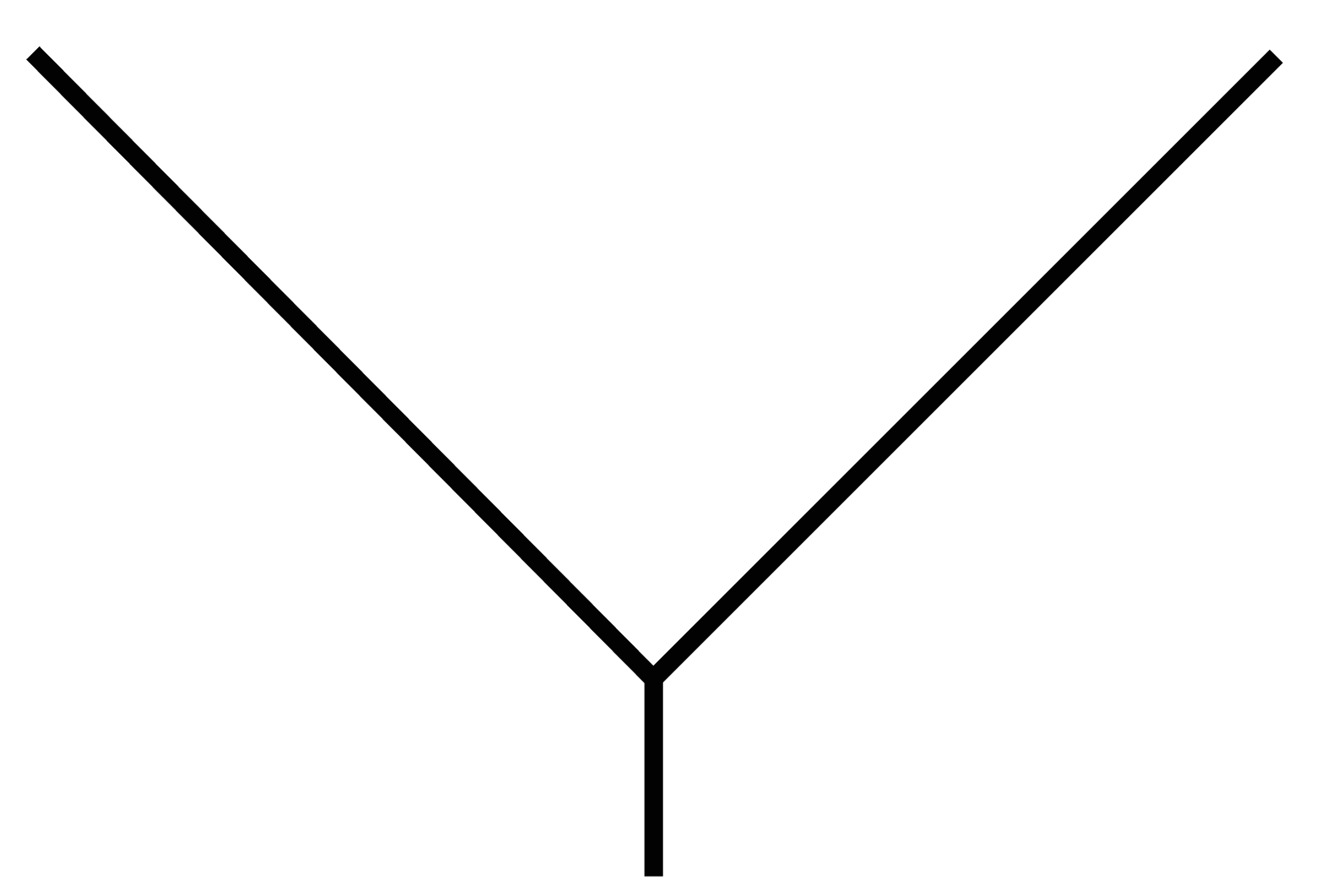} \right\}, 
\mathfrak{T}_2 = \left\{\arb{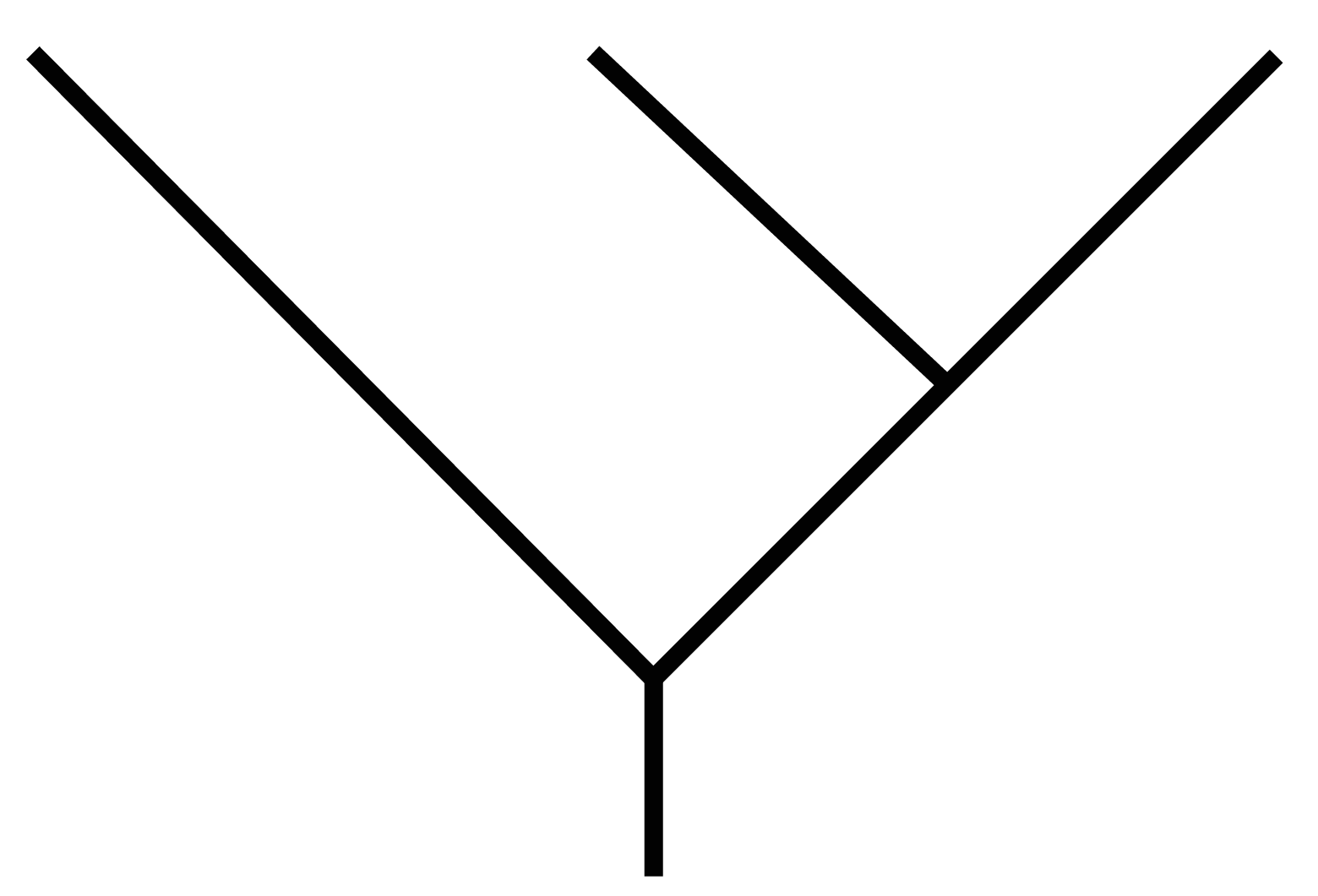}, \arb{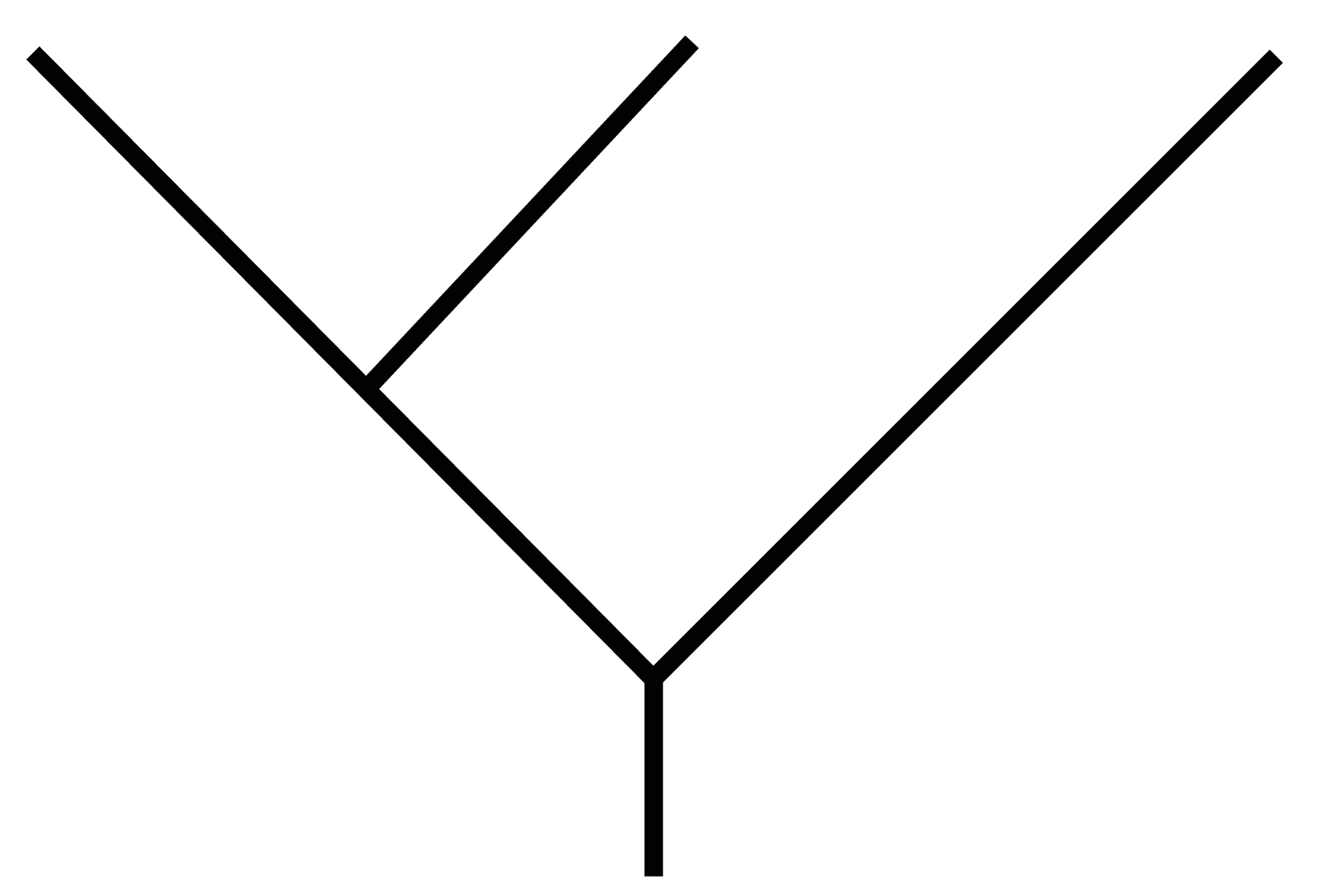}\right\},\\
\mathfrak{T}_3 &= \left\{ \arb{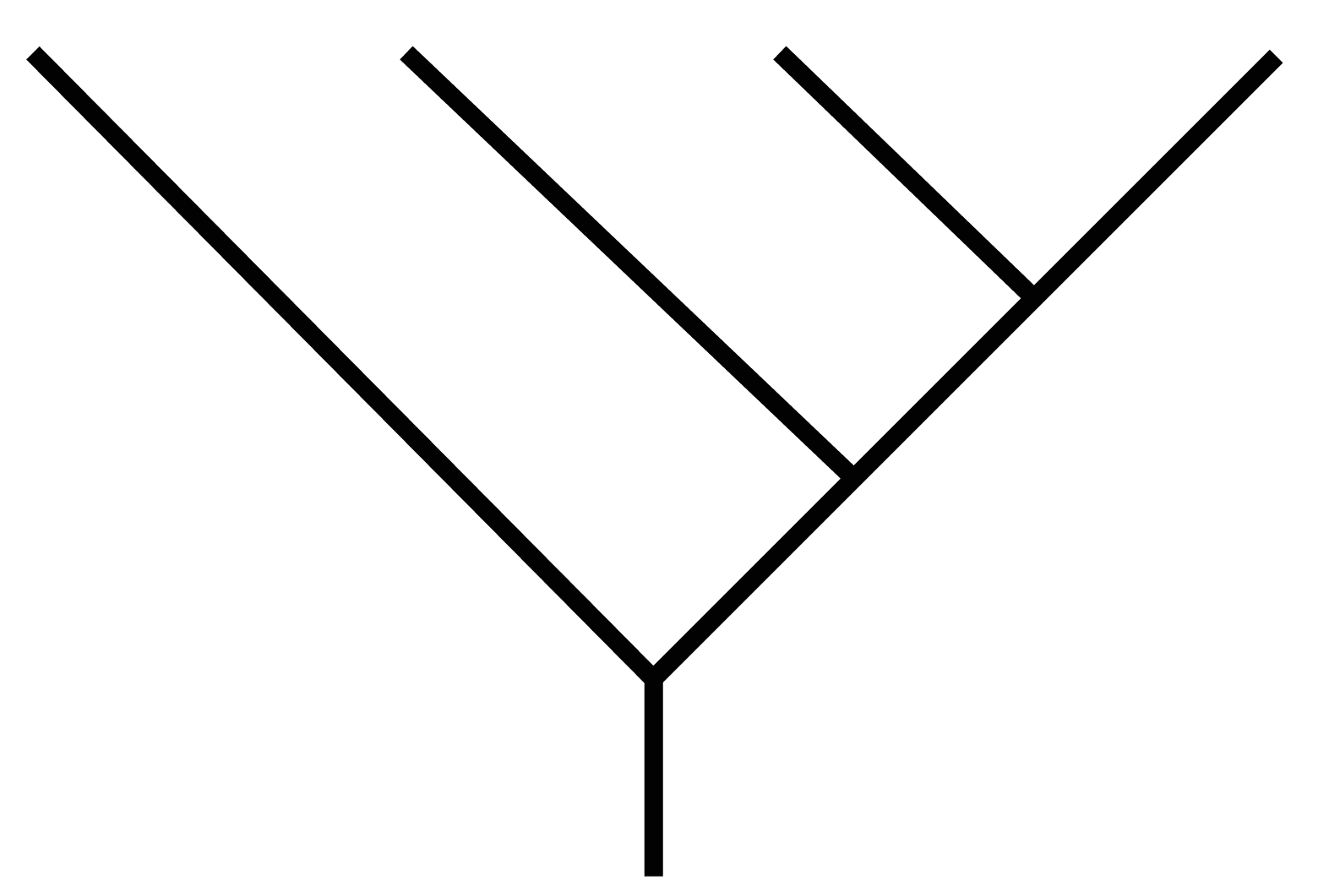}, \arb{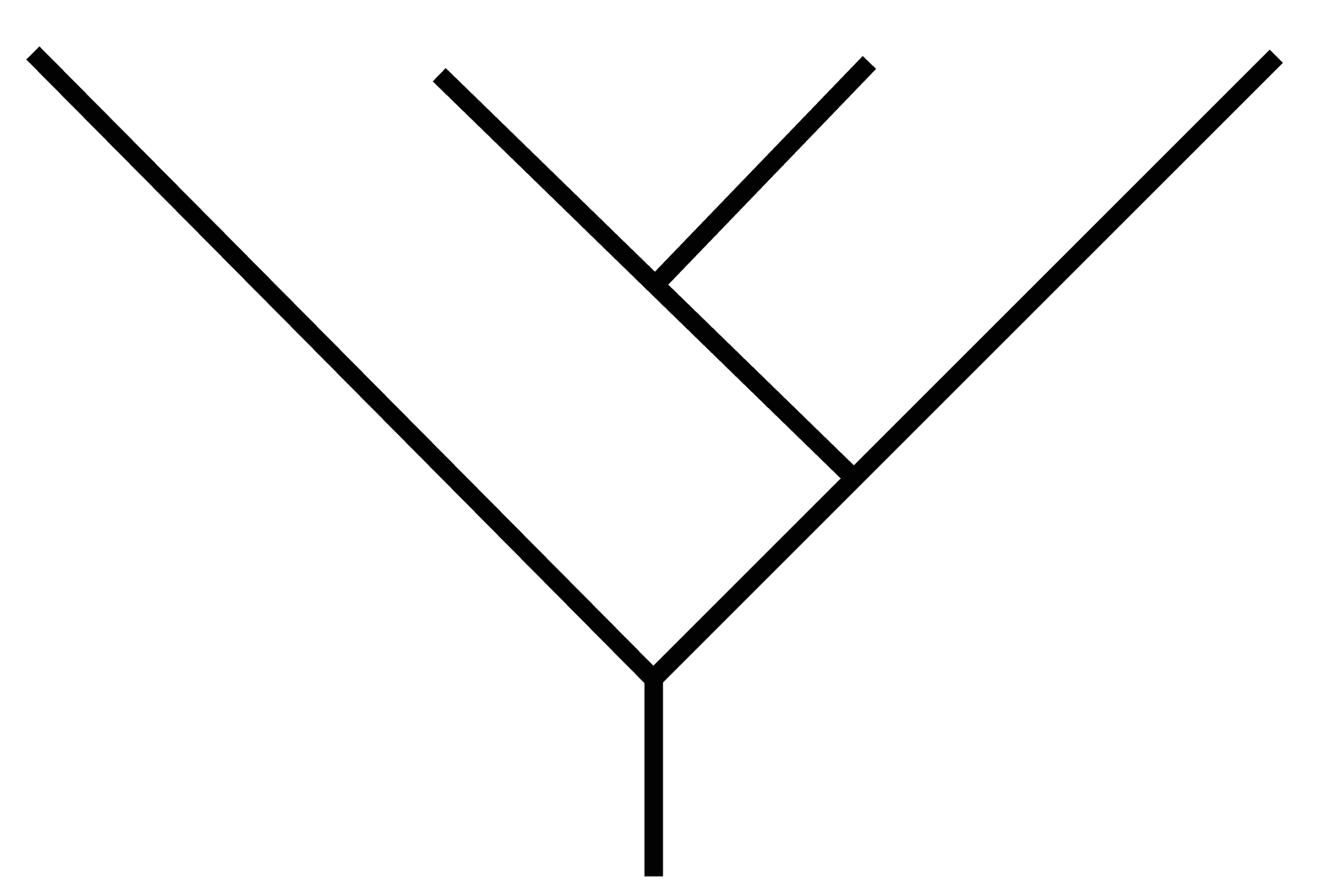}, \arb{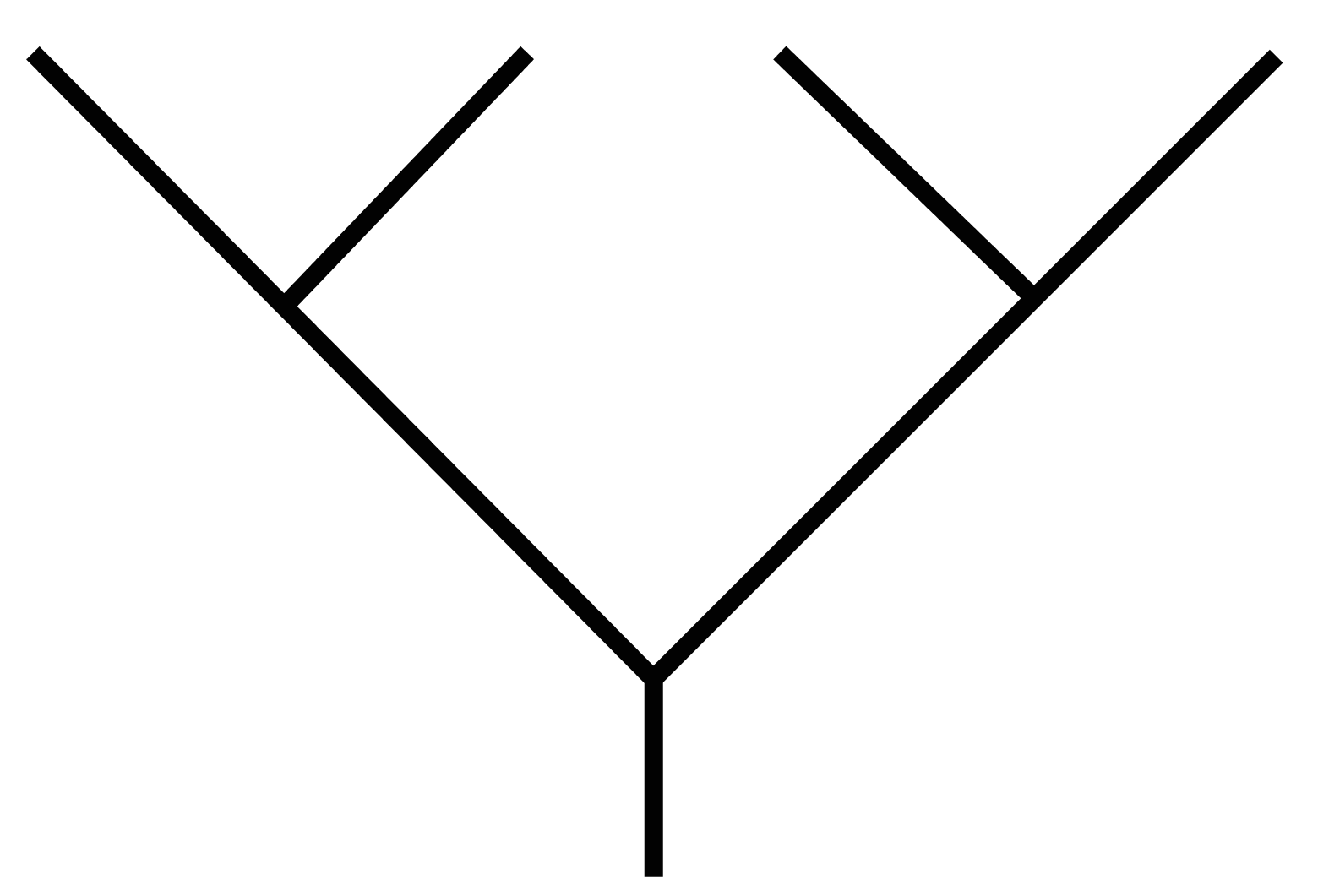}, \arb{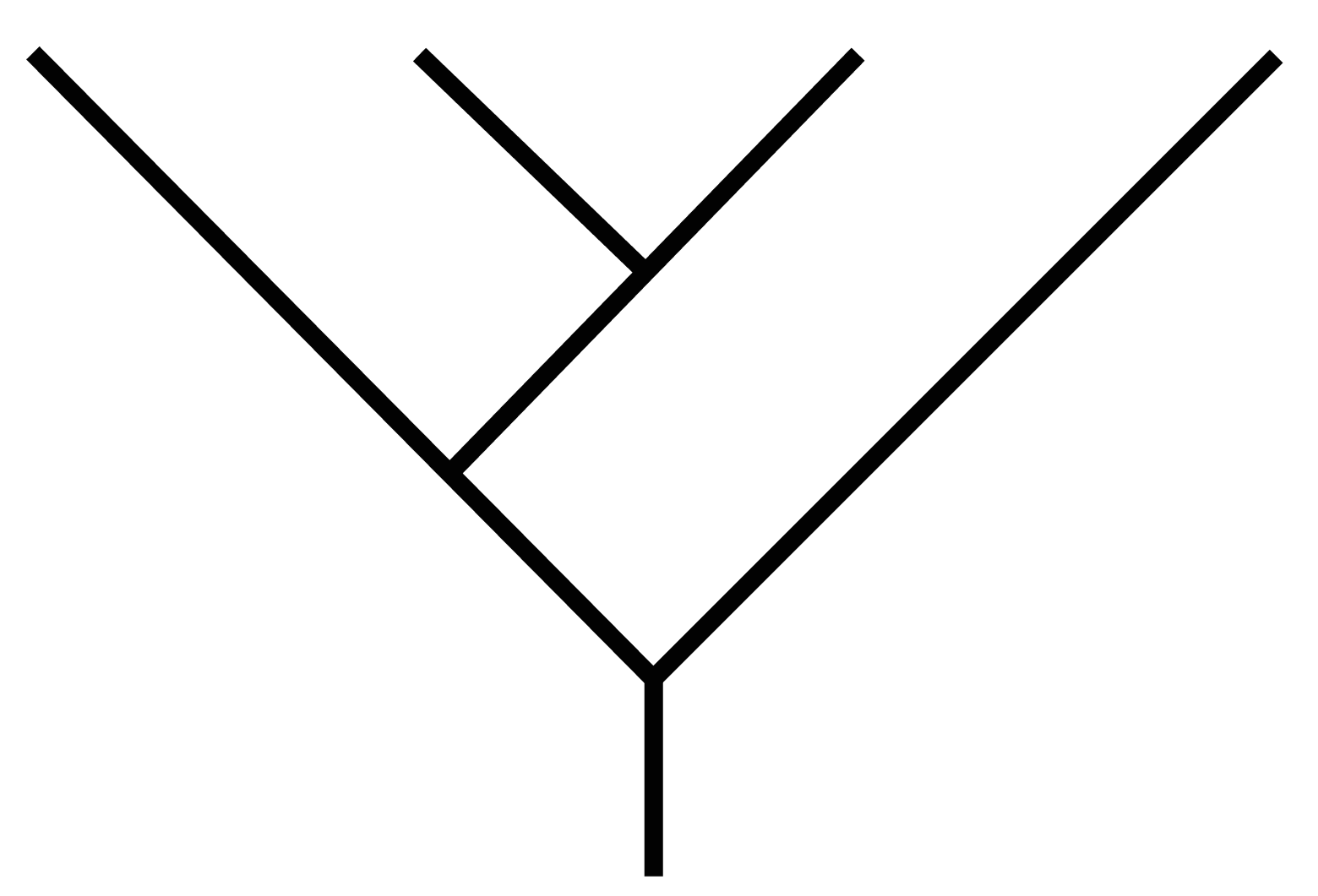}, \arb{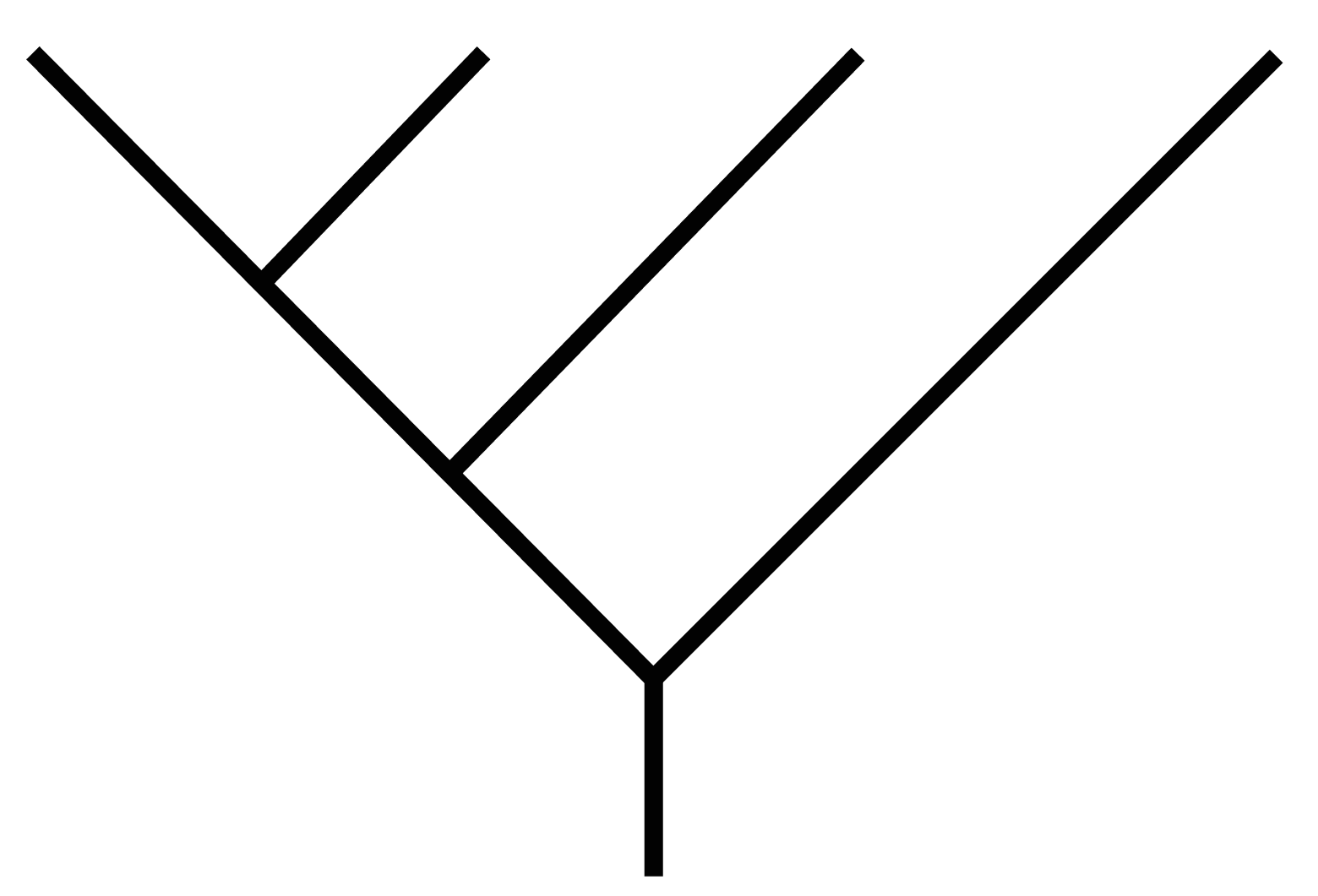} 
\right\}.
\end{align*}
The tree $|$ is the \emph{trivial tree}. A well known 
fact about planar binary trees is their cardinality 
$\#(\mathfrak{T}_n)=C_n:=\frac{1}{n+1}\binom{2n}{n}$, which is the 
$n$-th Catalan number. It is also known that the number of
ways of associating $n$ applications of a binary operator (e.g., balanced 
parenthesization) is $C_n$. Thus, if trees are suitably \emph{decorated} 
with a set of symbols, then there is a one-to-one correspondence between trees 
and dendriform words. Trees are decorated by attaching symbols to every interior 
vertex. 
\begde
Let $V_{int}$ be the set of interior vertices of tree $\tau\in \mathfrak{T}$ 
and $D$ a finite set of symbols. A decoration of $\tau$ is any injection $\rho: 
V_{int} \rightarrow D$.
\endde

\begin{figure}[H]
\caption{Tree decoration}
\label{fig:tree_labeling}
\begce
\includegraphics[width=2.5cm]{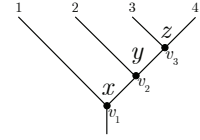}
\endce
\end{figure}

\begex Let $\tau=\arbsmall{321}$, $D=\{x,y,z\}$ and $V_{int}=\{v_1,v_2,v_3\}$, 
where $v_i$ is the vertex where the paths starting from leaves $i$ 
and $i+1$ join together. Figure \ref{fig:tree_labeling} shows the decoration of 
$\tau$ by $\rho$, where $\rho(v_1)=x$, $\rho(v_2)=y$ and $\rho(v_3)=z$.
\endex
\vspace*{-0.2in}
\vspace*{-0.1in}

In general, any $\tau\in \mathfrak{T}_n$ can be decorated by the letters in 
the word $\eta=x_{i_1}\cdots x_{i_n}\in X^n$, i.e., $\rho(v_j)=x_{i_j}$. The 
set of all trees decorated by $X^\ast$ is denoted as $\mathfrak{TD}X^\ast$, 
$\abs{\tau}$ is the order of $\tau\in \mathfrak{TD}X^\ast$, 
and the \emph{foliation} of $\tau$ is the mapping $\psi: \mathfrak{TD}X^\ast 
\rightarrow X^\ast$. A convenient way to consider the decoration of a tree with 
the letters of a word is by defining the operation $(\cdot;\cdot) : X^\ast 
\times \mathfrak{T} \rightarrow  \mathfrak{TD}X^\ast :(\eta;\tau) 
\mapsto\tau_\eta$. The notation $\tau_\eta$ makes explicit the fact that a tree 
$\tau \in \mathfrak{T}$ is being decorated by the word $\eta \in X^\ast$. For 
example, 
\begdi
(x_ix_jx_k;\arbsmall{321})=\arblarge{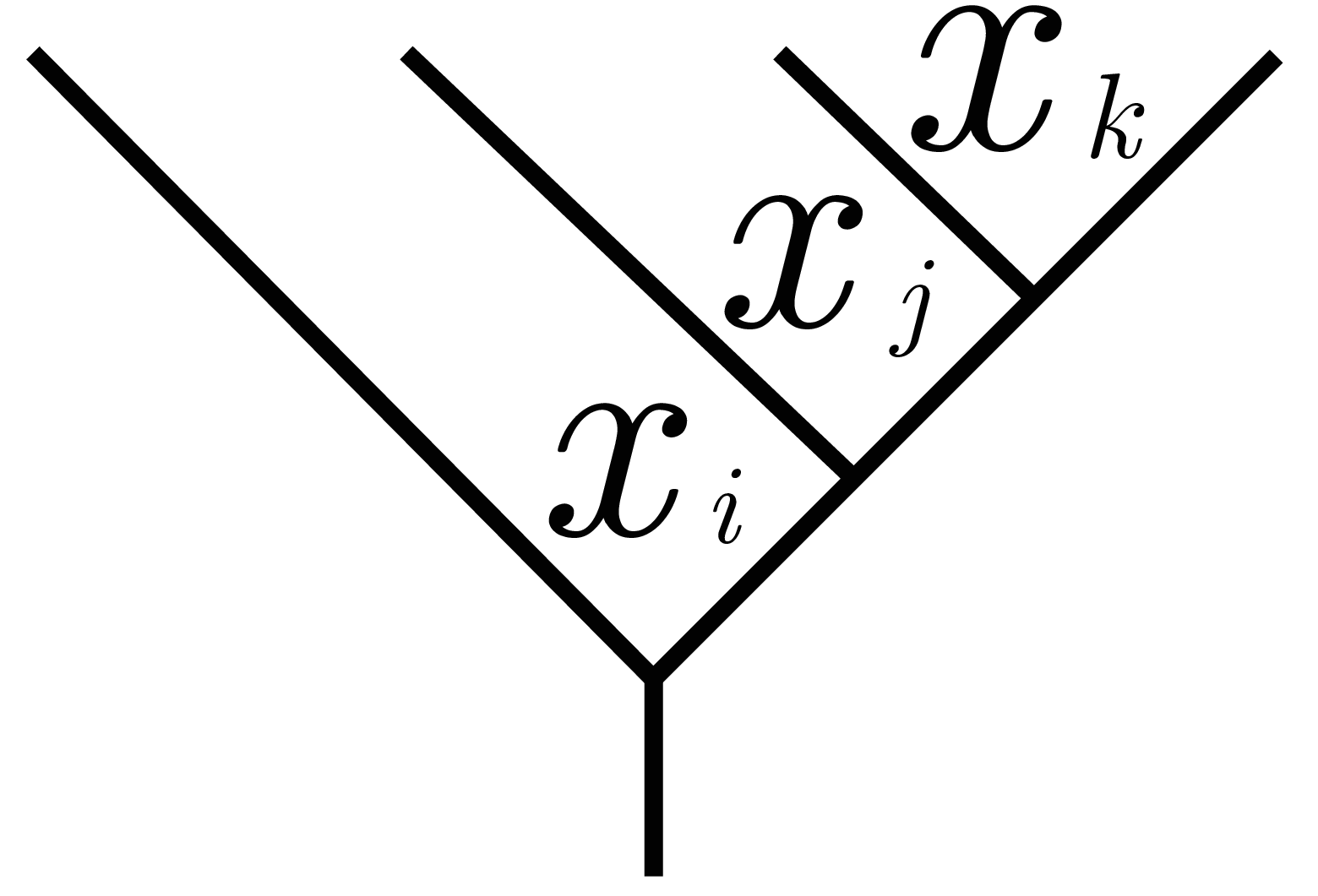}
\enddi
A formal power series on decorated trees is any mapping 
$c:\mathfrak{TD}X^\ast\rightarrow \re^{\ell\times n}\!:\eta\mapsto (c,\eta)$. 
The set of formal power series on decorated trees is $\allseriesTDXelln$, and 
forms an $\re$-vector space. The subset of series 
with finite support is $\allpolyTDXelln$. In this context, the decoration of 
trees is a bilinear operation. 

One way of constructing new trees from a given set of trees (usually 
called a \emph{forest}) is by the operation of \emph{grafting}.
\begde The \emph{grafting} of trees is an $n$-ary operation $\vee$ consisting 
of joining together $n$ trees to the same root to form a new tree. More 
precisely, $ \vee : \underbrace{\mathfrak{T} \times  \cdots \times 
\mathfrak{T}}_{\displaystyle \mbox{$n$ times}}  \rightarrow \mathfrak{T}$ such 
that
\begin{align*}
\vee(\tau^1 ,  \cdots  , \tau^n) =\overset{\;\;\;\;\tau^1\; \tau^2\; 
\tau^3\quad \;\; \tau^{n-1} \tau^n}{\arbsize{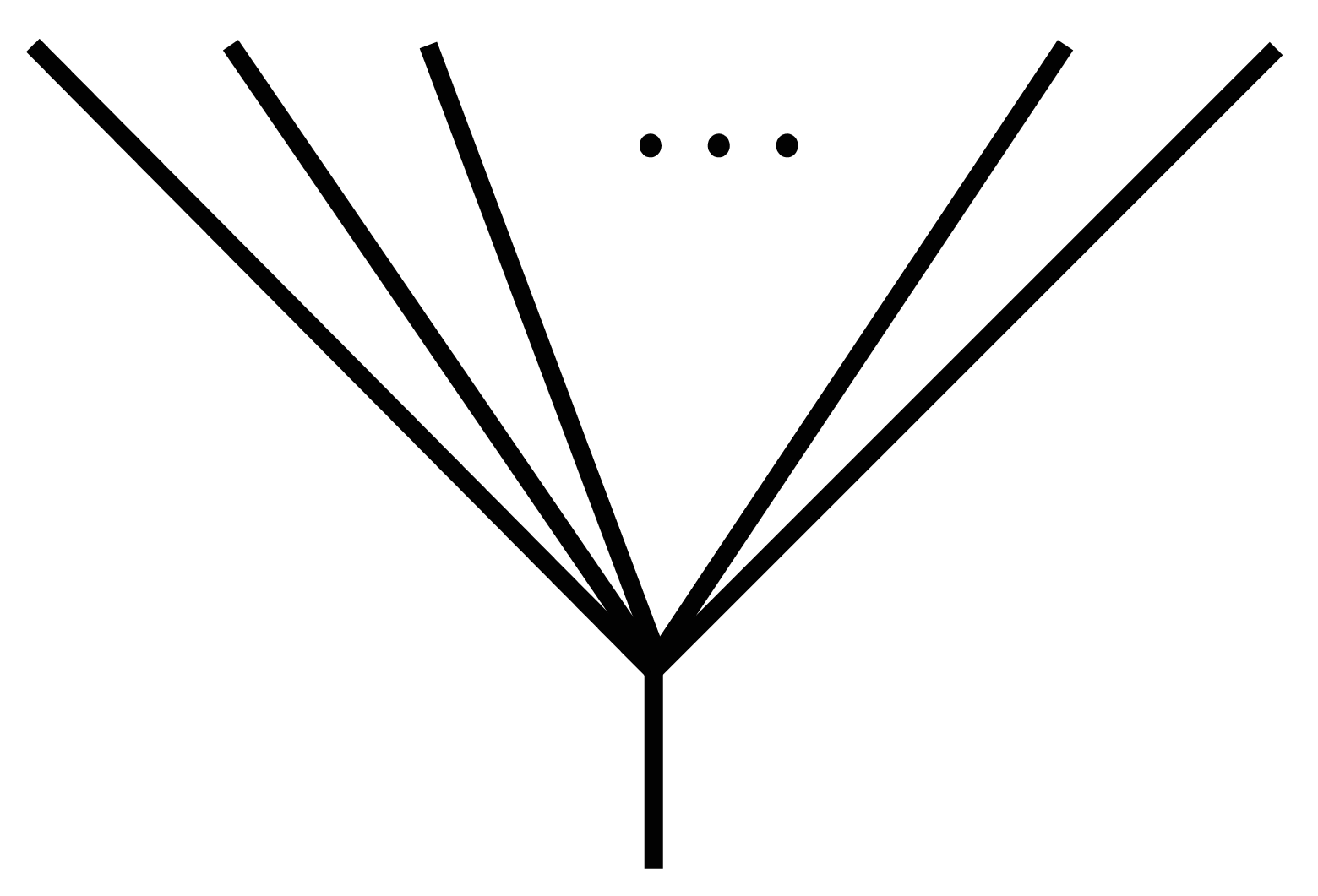}}.
\end{align*}
\endde

\begin{figure}[H]
\caption{$a)$ left comb, $b)$ right comb}
\label{fig:left-right-combs}
\begce
\includegraphics[width=5.cm]{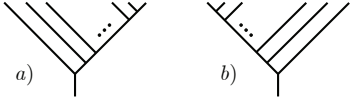}
\endce
\end{figure}

Grafting is for trees what catenation is for words. For example, 
\begin{align*}
& \arb{t} \vee \arb{t} = \arb{1}, \quad \arb{1} \vee \arb{1} = \arb{131}.
\end{align*}
Observe that if $\tau = \vee(\tau^1   \cdots   \tau^m)  \in\mathfrak{T}_n$, 
then $\tau^{i} \in\mathfrak{T}_{m_i}$ with $\sum_i m_i=n-1$. In this 
paper, the focus is on binary grafting, i.e., $m=2$. Any planar binary tree 
can be decomposed uniquely as $\tau=\tau^1\vee \tau^2$ since by definition 
any planar binary tree interior vertex is trivalent (has one root and two 
leaves). Other tree decompositions such as the ones used in Hopf algebras of 
trees are not unique \cite{Holtkamp_2011}. The tree $\tau^1$ (respectively 
$\tau^2$) is the left part (respectively the right part) of $\tau$. Further 
decompositions allow one to write any planar binary tree in terms of the 
trivial tree $|$. The grafting operation $\vee$ makes $\mathfrak{T}$ the free 
magma algebra with one generator. It is neither commutative nor associative but 
is of order one with respect to the grading in terms of internal vertices. That 
is, for two trees $\tau^1$, $\tau^2$ of order $n_1$, $n_2$, respectively, the 
product $\tau^1 \vee \tau^2$ is of order $n_1 + n_2 + 1$. Particular types of 
trees that allow an easy decomposition are the so-called right-comb and 
left-comb as shown in Figure \ref{fig:left-right-combs}.
\vspace*{-0.1in}
\vspace*{-0.1in}
\noindent Clearly for a right-comb (respectively left-comb) 
$\tau_r^n=\tau_r^{n-1}\vee |$ (respectively $\tau_l^{n}=|\vee \tau_l^{n-1}$), 
where $\tau_r^k$ denotes the $k$-th order right-comb (respectively 
$\tau_l^k$ denotes the $k$-th order left-comb).

One way of realizing the decoration of a tree is by attaching a 
letter from the alphabet $X$ to every grafting operation used in the 
construction, say $\vee_{x_i}$. For example,
\begin{subequations} \label{subeq:non_commuting_grafting}
\begin{align} 
&\arb{t} \vee_{x_i} \arb{t} = \arblarge{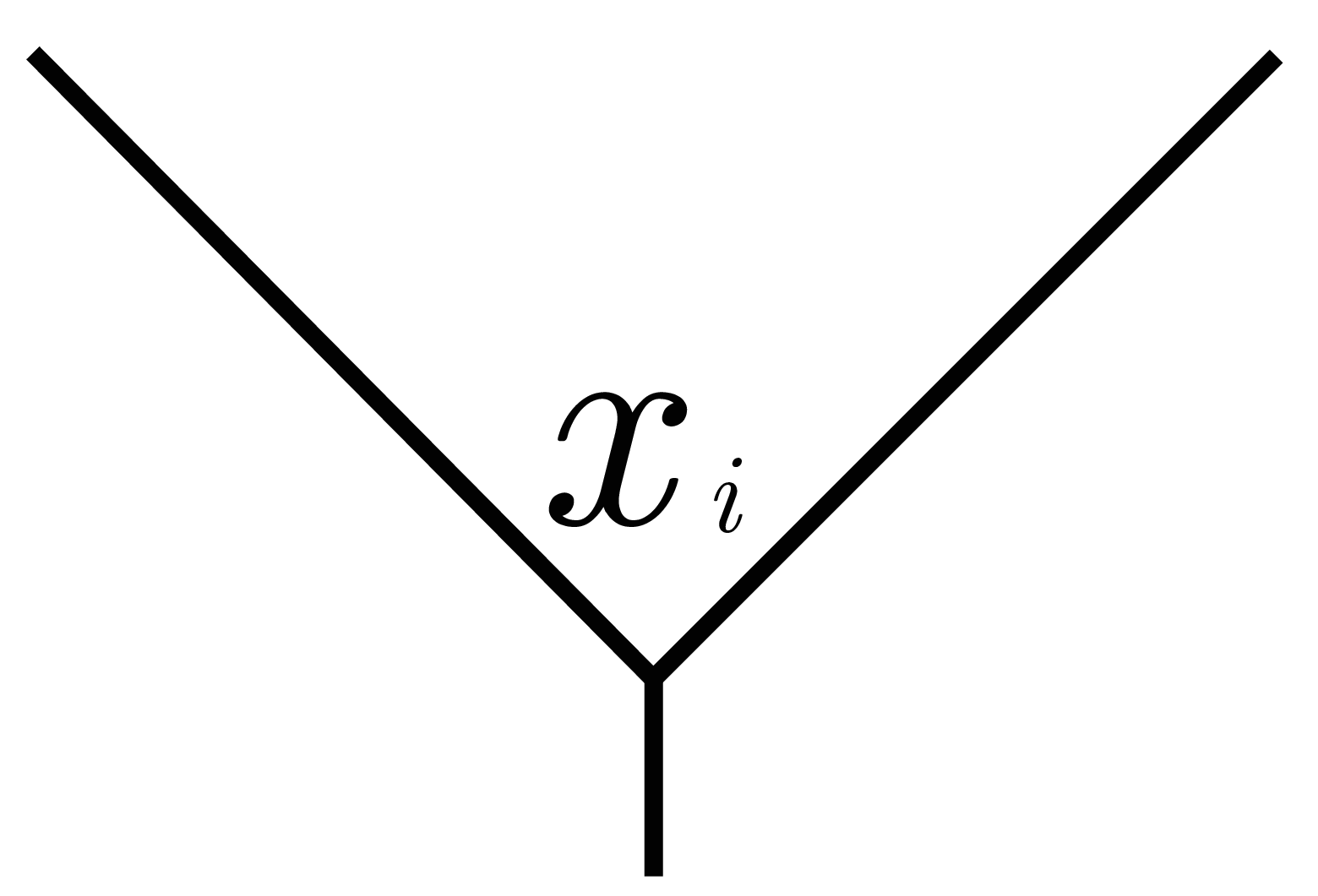} \\
&\arb{t} \vee_{x_i} \left( \arb{t} \vee_{x_j} \arb{t} \right)  = \arb{t} 
\vee_{x_i} \arblarge{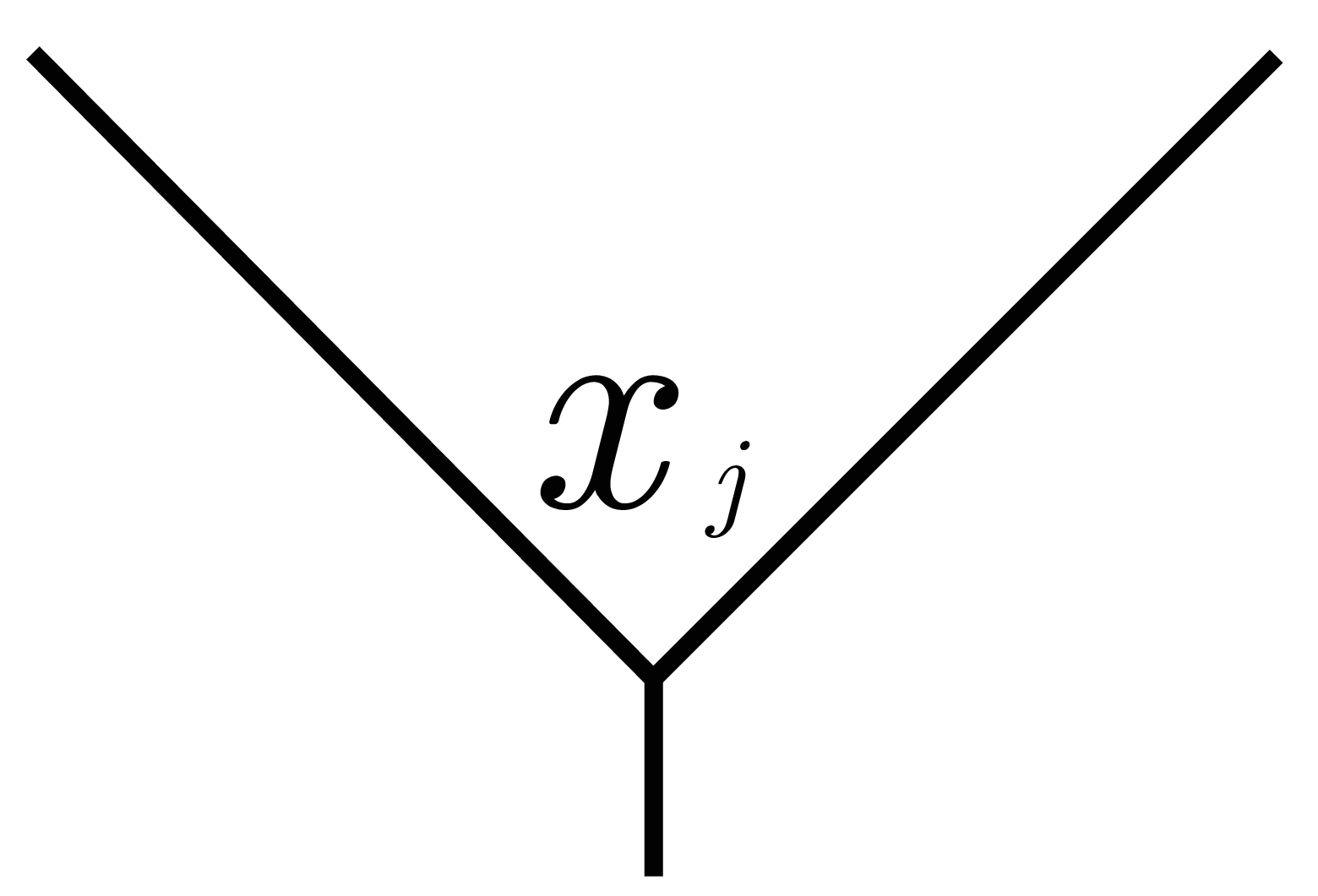} =  \arblarge{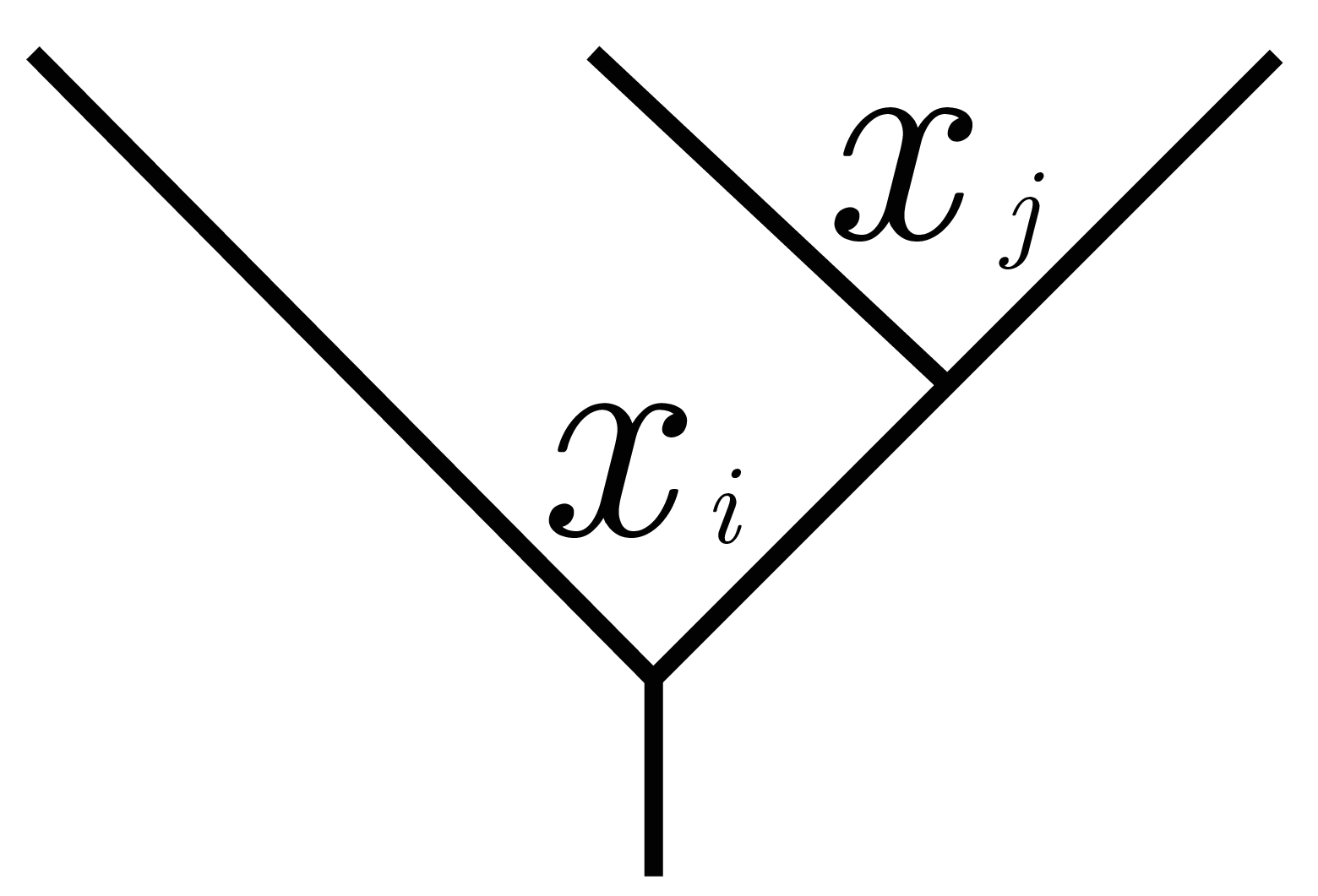} \\
 &\left( \arb{t} \vee_{x_i} \arb{t} \right) \vee_{x_j} \left( \arb{t} 
\vee_{x_k} \arb{t} \right)   = \arblarge{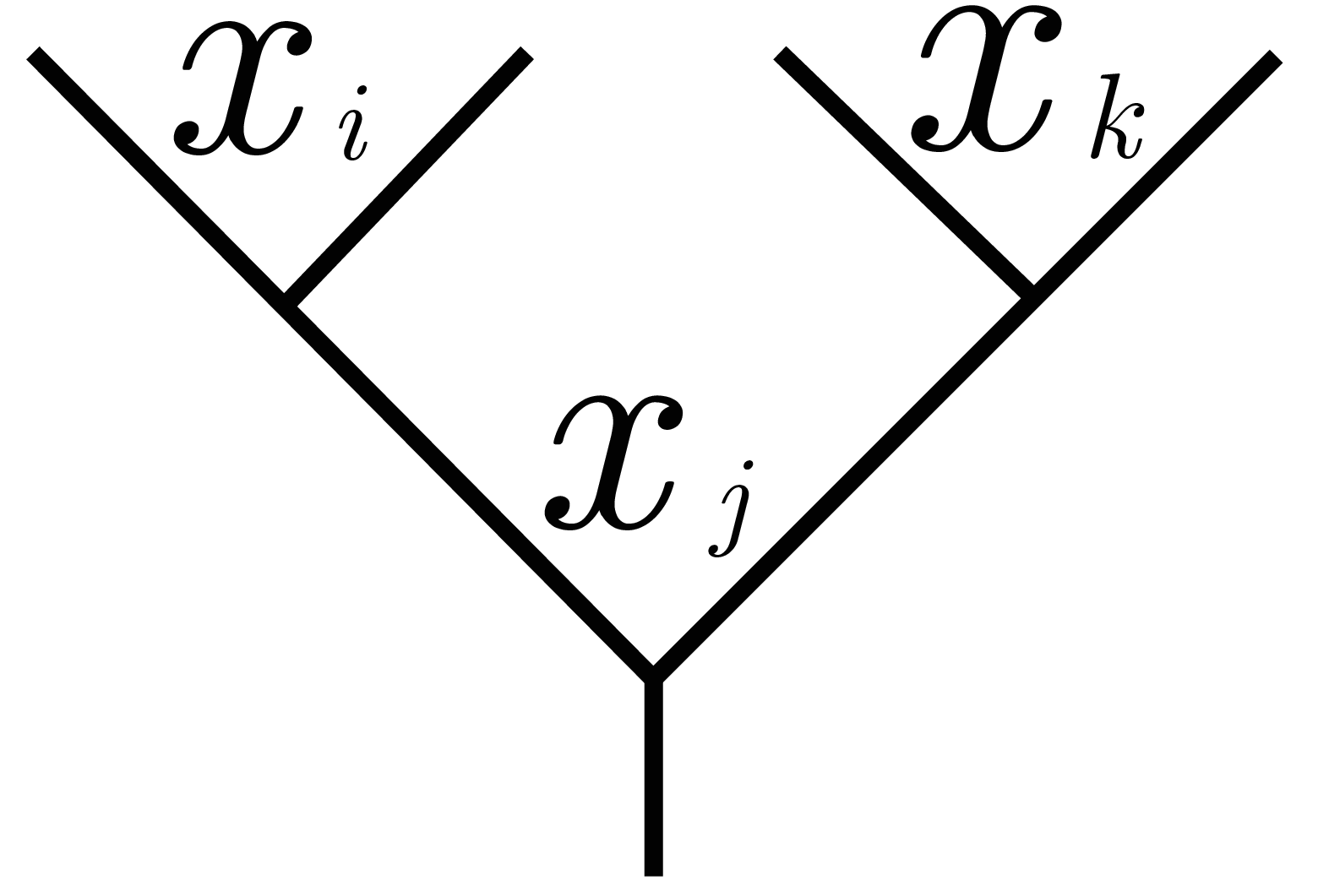}.
\end{align}
\end{subequations}

The grafting operation allows an explicit description of the correspondence 
between the sets $\mathfrak{T}X^\ast$ and $\mathfrak{TD}X^\ast$. This is 
provided by the isomorphism $\Phi:\mathfrak{T}X^\ast 
\rightarrow \mathfrak{TD}X^\ast$ with the inductive definition 
\begin{align*}
\Phi(\eta_{\tau}) = \left\{\begin{array}{ll}
               |\; \vee_{x_i}  \Phi(\eta'_{\tau'}), & \mbox{if} \quad \eta_\tau 
= x_i \prec \eta'_{\tau'}, \\
		\Phi(\eta'_{\tau'}) \vee_{x_i} | , & \mbox{if} \quad \eta_\tau 
= \eta'_{\tau'} \succ x_i, \\
\Phi(\eta'_{\tau'}) \vee_{x_i} \Phi(\eta''_{\tau''}) , & \mbox{if} \quad 
\eta_\tau = 
\eta'_{\tau'} \succ x_i \prec \eta''_{\tau''}, 
               \end{array}\right.
\end{align*}
where $x_i\in X$, $\eta'_{\tau'},\eta''_{\tau''}\in \mathfrak{T}X^\ast$, 
$\Phi(\emptyset)=|$ and $\Phi(x_j)=| \;\vee_{x_j} |$ for $x_j\in X$. 
For example,
\begin{align*}
&\Phi\left(\emptyset\right) = \arblarge{t},  \Phi\left(x_i \right) = 
\arblarge{1x}, \Phi\left(x_i\prec x_j \right) = \arblarge{21xy}, \\
& \Phi\left(x_i\succ (x_j \prec x_k) \right) = \arblarge{131xyz}.
\end{align*}

In \cite{Melancon_92}, the free magma $\mathfrak{T}X^\ast$ is defined directly 
as the set of all planar binary trees whose leaves are decorated with the 
letters in $X$. Any $\eta\in\mathfrak{T}X^\ast$ will be denoted as 
$\eta_{\tau}$, where it is made explicit the fact that for any dendriform word 
there exist a decorated tree ${\tau}\in\mathfrak{TD}X^\ast$ providing the order 
in which the products $\prec$ and $\succ$ appear. The corresponding tree is 
then obtained as $\Phi(\eta_\tau)=\tau_\eta \in \mathfrak{TD}X^\ast$, and its 
inverse satisfies $\Phi^{-1}(\tau_\eta)=\eta_\tau \in \mathfrak{T}X^\ast$. 
Moreover, the foliation of $\tau_\eta$ can be written in terms of the map 
$\varphi$ of dendriform words as $\psi(\tau_\eta) = 
\varphi(\Phi^{-1}(\tau_\eta))=\eta \in X^\ast$. The isomorphism $\Phi$ 
is extended linearly in the natural way over $\allseriesTXelln$. Hereafter, due 
to the isomorphism between $\mathfrak{T}X^\ast$ and $\mathfrak{TD}X^\ast$, no 
notational distinction will be made between products in $\mathfrak{T}X^\ast$ and 
products in $\mathfrak{TD}X^\ast$.

The grafting operation allows one to define the \emph{shuffle product of trees} 
in an inductive manner. 

\begde \cite{Loday-Ronco_98} Let $\tau^1= \tau^{11}\vee_{x_i} \tau^{12}$ and 
$\tau^2=\tau^{21}\vee_{x_j} \tau^{22}$ with $\tau^{k1},\tau^{k2} \in 
\mathfrak{TD}X^\ast$. The recursive definition of 
$\shuffleNC:\mathfrak{TD}X^\ast\times \mathfrak{TD}X^\ast \rightarrow 
\allpolyTDX $ is 
\begeq \label{eq:shuffleNC}
\tau^1 \shuffleNC \tau^2 = \tau^{11} \vee_{x_i} ( \tau^{12} \shuffleNC \tau^2 ) 
+ (\tau^1 \shuffleNC \tau^{21}) \vee_{x_j} \tau^{22},
\endeq
where $| \shuffleNC \tau = \tau \shuffleNC | = \tau$ for any $\tau\in 
\mathfrak{TD}X^\ast$.
\endde

\begex \label{ex:non-commutative_shuffle} The shuffle product 
\begin{align*}
\lefteqn{\arb{21xy} \shuffleNC \arb{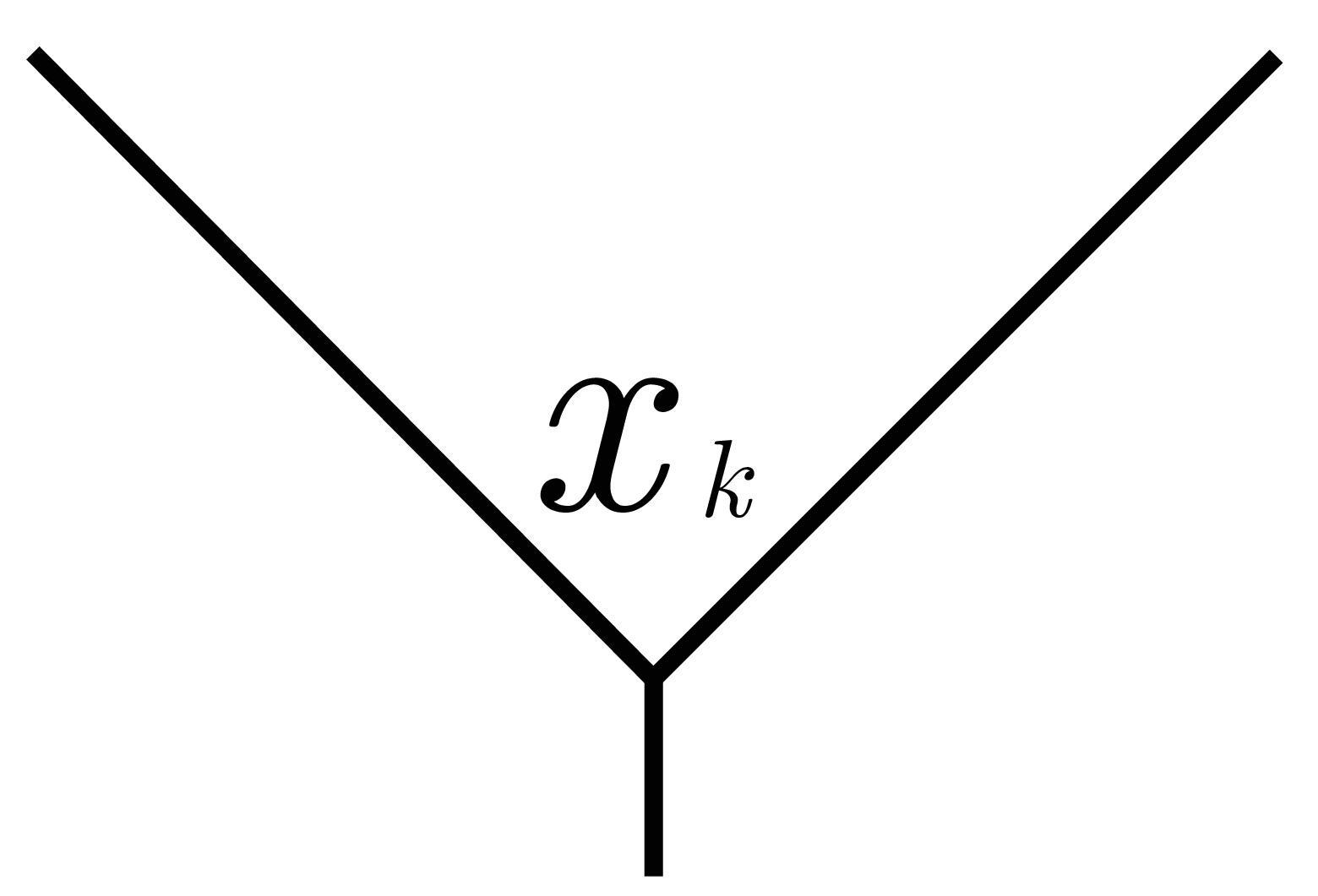} =  \left( \arb{t} \vee_{x_i} \arb{1y} 
\right)   \shuffleNC \left( \arb{t} \vee_{x_k} \arb{t} \right)}\\
 & = \arb{t} \vee_{x_i} \left( \arb{1y} \shuffleNC \arb{1z} \right) + 
 \left( \arb{21xy} \shuffleNC \arb{t} \right) \vee_{x_k} \arb{t}\\
 & = \arb{t} \vee_{x_i} \left( \arb{t} \vee_{x_j}  \arb{1z} + \arb{1y}  
 \vee_{x_k} \arb{t}  \right) + \arb{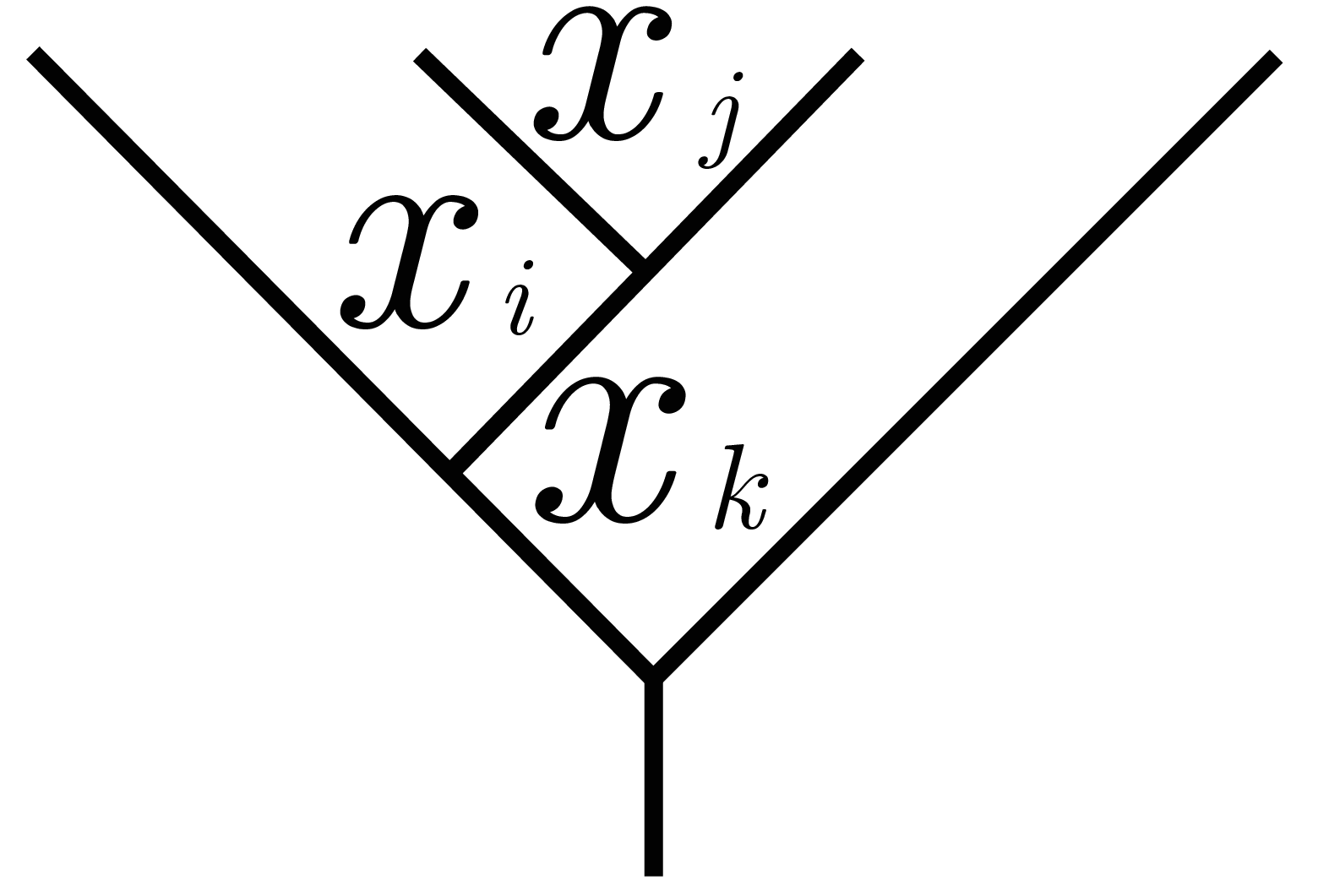} \\
 & = \arb{321xyz} + \arb{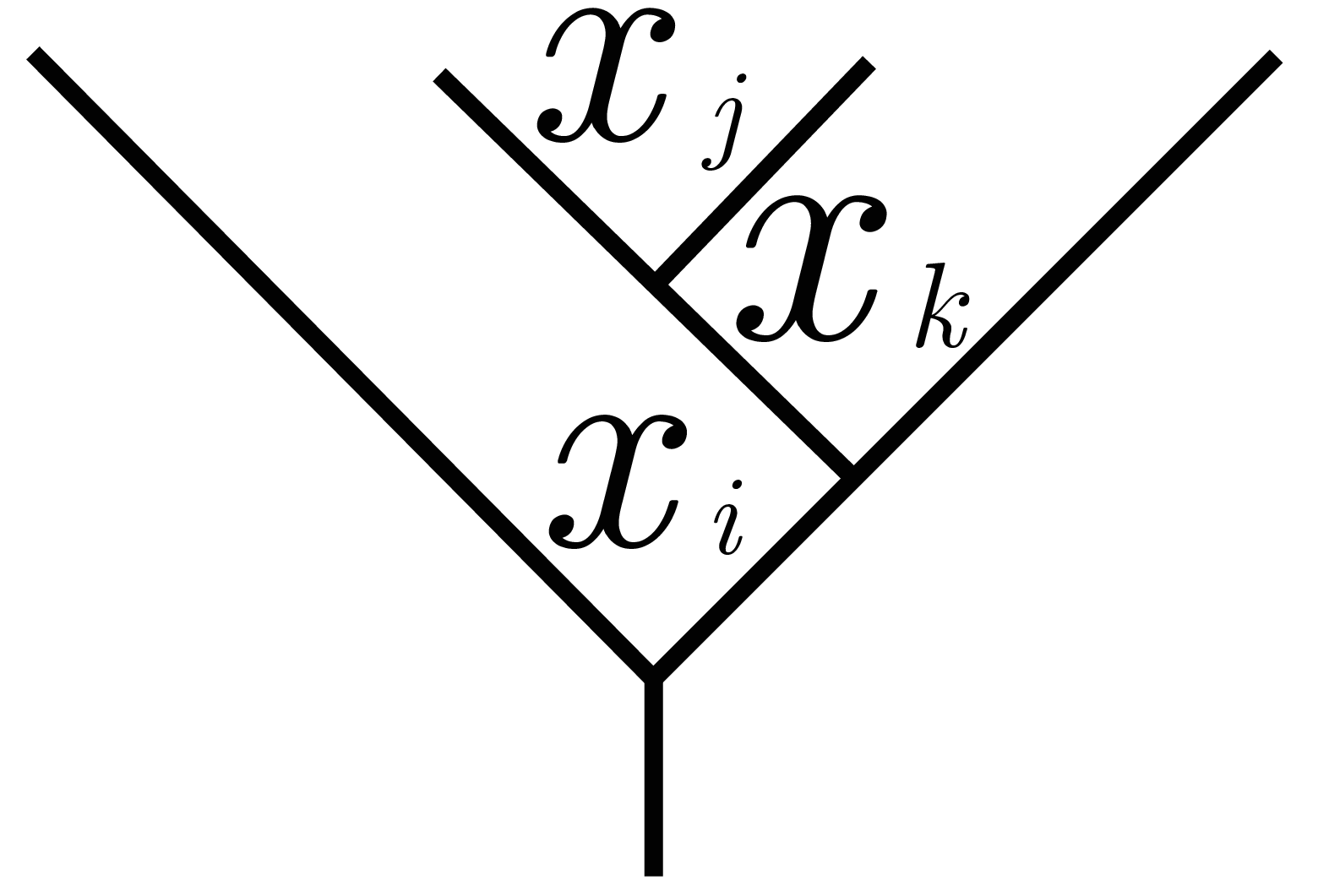} + \arb{213xyz}.
\end{align*}
\endex

\begle \cite{Loday-Ronco_98}
The shuffle product $\shuffleNC$ of trees is non-commutative and associative.
\endle

The dendriform products for trees are given next. 

\begde \cite{Loday-Ronco_98} Let $\tau^1= \tau^{11}\vee_{x_i} \tau^{12}$ and 
$\tau^2= \tau^{21}\vee_{x_j} \tau^{22}$. The dendriform products for 
decorated trees are: 
\begin{subequations} \label{subeq:dendriform_tree}
\begin{align}
\label{eq:dendriform_tree1} \tau^1 \prec \tau^2 & = 
\tau^{11}\vee_{x_i} 
(\tau^{12}\shuffleNC \tau^2) \\
\label{eq:dendriform_tree2} \tau^1 \succ \tau^2 & = 
(\tau^{1}\shuffleNC \tau^{21}) \vee_{x_j} \tau^{22}. 
\end{align}
\end{subequations}
\endde
\vspace*{0.05in}
\begth \cite{Loday-Ronco_98} \label{th:three_dendriform_products}
The products $\prec$ and $\succ$ in \rref{subeq:dendriform_tree} satisfy the 
axioms of dendriform algebras given in \rref{subeq:dendriform_axioms}.
\endth

From Theorem \ref{th:three_dendriform_products}, \rref{eq:shuffleNC} and 
\rref{subeq:dendriform_tree}, it is clear that $\prec + \succ = \shuffleNC $ as 
in Section \ref{sec:2}, and therefore, with the help of the mapping $\Phi$, the 
relationship between the dendriform and shuffle products on $\mathfrak{T}X^\ast$ 
and $\mathfrak{TD}X^\ast$ is:
\begin{align*}
\Phi(\eta_{\tau^1} \prec \xi_{\tau^2}) & = \Phi(\eta_{\tau^1}) \prec 
\Phi(\xi_{\tau^2})=\tau^1 \prec \tau^2 \\
\Phi(\eta_{\tau^1} \succ \xi_{\tau^2}) & = \Phi(\eta_{\tau^1}) \succ 
\Phi(\xi_{\tau^2})= \tau^1 \succ \tau^2\\
\Phi(\eta_{\tau^1} \shuffleNC \xi_{\tau^2}) & = \Phi(\eta_{\tau^1}) \shuffleNC 
\Phi(\xi_{\tau^2})= \tau^1 \shuffleNC \tau^2\\
\Phi(\eta_{\tau^1} \shuffleNC \emptyset) & = \Phi(\emptyset \shuffleNC 
\eta_{\tau^1}) = \tau^1 \shuffleNC |.
\end{align*}

This subsection ends with two key lemmas employed in Section \ref{subsec:3-B} to 
characterize the grouping of non-commutative iterated integrals.

\begle \label{le:characteristic_series_tree_decomp} For any $n\ge 0$,
\begeq \label{eq:characteristic_series_tree_decomp}
\hspace*{-0.in}\mbox{char}(\mathfrak{T}_{n+1}) := \sum_{\tau \in 
\mathfrak{T}_{n+1}} \!\!\!\tau = \sum_{i=0}^{n} \mbox{char}(\mathfrak{T}_{n-i}) 
\vee 
\mbox{char}(\mathfrak{T}_i).
\endeq
\endle
\vspace*{0.05in}

\begpr
Recall that $\mathfrak{T}_{n}$ is 
$\#\mathfrak{T}_{n}=C_n$. Since the grafting operation is non-commutative and 
provides a unique decomposition of planar binary trees, one can prove the claim 
by showing that the right-hand side of 
\rref{eq:characteristic_series_tree_decomp} produces a number of summands equal 
to the $(n+1)$ Catalan number. First, note that the grafting operation does 
not generate extra trees in the sense that 
\begin{align*}
\# \supp\left(\sum_{i=1}^{n_1} \tau_{1,i} \vee 
\sum_{j=1}^{n_2}\tau_{2,j}\right) & = n_1 n_2.
\end{align*}
It then follows that
\begin{align*}
\#\supp\left(\mbox{char}(\mathfrak{T}_{n+1})\right) & = \sum_{i=0}^{n} 
\#\supp\left(\mbox{char}(\mathfrak{T}_{n-i}) \vee 
\mbox{char}(\mathfrak{T}_i)\right) \\
& =\sum_{i=0}^{n} C_{n-i}C_{i} = C_{n+1},
\end{align*}
which is Segner's recurrence relation for the $(n+1)$ Catalan number 
\cite{Segner_1759}.
\endpr

The collection of all trees of a certain order can be described in 
terms of the non-commutative shuffle product.
\begle \label{le:characteristic_shuffle_identity} The summation of all 
undecorated trees of order $n\ge 0$ is given by
\begeq \label{eq:characteristic_shuffle_identity}
\mbox{char}(\mathfrak{T}_n) = \arb{1}^{\shuffleNC n},
\endeq
where $\arbsmall{1}^{\shuffleNC n+1} = (\arbsmall{1}^{\shuffleNC 
 n})\shuffleNC \arbsmall{1}$ and $\arbsmall{1}^{\shuffleNC 0}=\emptyset$.
\endle

\begpr The proof is done by induction on the number of shuffles. For $n=0,1$, 
the identity holds trivially. For $n=2$, it is easy to see that
\begdi
\arb{1} \shuffleNC \arb{1} = \arb{12}+\arb{21} = \mbox{char}(\mathfrak{T}_2).
\enddi
Assume now that \rref{eq:characteristic_shuffle_identity} holds up to some 
$n\ge 1 $. Using Lemma~\ref{le:characteristic_series_tree_decomp} and the 
associativity of $\shuffleNC$, it follows that
\begin{align*}
\lefteqn{\hspace*{-0.3in}\arb{1}^{\shuffleNC (n+1)} = \left(\arb{1}^{\shuffleNC 
n}\right) 
\shuffleNC \arb{1}} \\
= & \sum_{i=0}^{n-1} \mbox{char}(\mathfrak{T}_{i}) \vee 
\underbrace{\left(\mbox{char}(\mathfrak{T}_{n-1-i})\shuffleNC 
{\arb{1}}\right)}_{\displaystyle \mbox{char}(\mathfrak{T}_{n-i})}   \\
  & + \left(\left(\sum_{i=0}^{n-1} \mbox{char}(\mathfrak{T}_{i}) 
\vee \mbox{char}(\mathfrak{T}_{n-1-i})\right) \shuffleNC \arb{t} \right) \vee 
\arb{t}.
\end{align*}
Given that $\tau\shuffleNC |=\tau$, and using the induction hypothesis, the 
last summand above is  
\begdi
\left(\sum_{i=0}^{n-1} \mbox{char}(\mathfrak{T}_{i}) 
\vee \mbox{char}(\mathfrak{T}_{n-1-i})\right) \vee \arb{t} = 
\mbox{char}(\mathfrak{T}_{n}) \vee \mbox{char}(\mathfrak{T}_{0}).
\enddi
Thus,
\begin{align*}
\lefteqn{\arb{1}^{\shuffleNC (n+1)}} \\
& = \sum_{i=0}^{n-1} \mbox{char}(\mathfrak{T}_{i}) 
\vee \mbox{char}(\mathfrak{T}_{n-i}) + \mbox{char}(\mathfrak{T}_{n}) 
\vee\mbox{char}(\mathfrak{T}_{0}) \\
& = \sum_{i=0}^{n} \mbox{char}(\mathfrak{T}_{i})\vee 
\mbox{char}(\mathfrak{T}_{n-i}) = \mbox{char}(\mathfrak{T}_{n+1}).
\end{align*}
\endpr

\subsection{Non-commutative iterated integrals} \label{subsec:3-B}

For a matrix-valued measurable function $u : [0,T] \rightarrow \re^{n\times 
q}$, define $\norm{u}_{L_1} = \int_{0}^{T} 
\norm{u(s)}_1 \,ds$. Note that $\norm{u(s)}_1=\max_j \{\sum_{i} 
\abs{u(s)_{ij}}\}$. Now let instead $u=(u_1,\hdots,u_m)$, where each $u_i:[0, 
T] \rightarrow \re^{n\times q}$. The norm for this $u$ is $\norm{u}:= 
\max_i \norm{u_i}_{L_1} $. The set $L^{m\times(n\times q)}_1[0,T]$ contains 
all measurable functions defined on $[0,T]$ having finite $\norm{\cdot}$ 
norm, and $B^{m\times (n\times q)}_1(R)[0,T]:= \{u\in L^{m\times 
(n\times q)}_1[0,T], \norm{u}\le R \}$. 

\begde
Let $u\in B^{m\times (n\times n)}_1(R)[0,T]$. The non-commutative iterated 
integral corresponding to $\eta_\tau \in \mathfrak{T}X^\ast$ for $t\in [0,T]$ 
is defined inductively by $E_{\emptyset}[u]=I$, and 
\begin{equation*}
E_{\eta_\tau}[u](t) = \int_{0}^t E_{{\xi}_{\tau^1}}[u](s) u_i(s) 
E_{{\nu}_{\tau^2}}[u](s)\, ds,
\end{equation*}
where $x_i\in X$, $\eta_\tau={\xi}_{\tau^1}\vee_{x_i}{\nu}_{\tau^2}$ with 
${\xi}_{\tau^1},{\nu}_{\tau^2}\in \mathfrak{T}X^\ast$, 
$\tau^1, \tau^2 \in \mathfrak{T}$, $u_0=I$, and $I$ denotes the identity 
matrix.
\endde

The mapping $E_{\eta_{\tau}}$ is extended linearly on $\allpolyTXnn$ in the 
natural way. For example, the iterated integrals corresponding to 
\rref{subeq:non_commuting_grafting} are, respectively, 
\begin{align*}
& E_{x_i}[u](t)=\int_0^t u_i(s)\,ds, \\
& E_{x_i \prec x_j}[u](t) = \int_0^t u_i(s) \int^s_0 u_j(\tau)\,d\tau ds, \\
& E_{x_i\succ (x_j \prec x_k)}[u](t)= E_{( x_i \succ x_j ) \prec x_k}[u](t)\\ 
& = \int_0^t  \left( \int_0^s u_i(\tau) 
\,d \tau \right)  u_j(s) \left( \int_0^s u_k(\tau)\,d\tau \right) ds. 
\end{align*}

For a planar binary tree $\tau=\tau^1 \vee \tau^2 \in \mathfrak{T}$, the 
\emph{tree factorial} is defined as 
\begdi
\gamma(\tau) = (\abs{\tau^1}+\abs{\tau^2}+1)\gamma(\tau^1)\gamma(\tau^2),
\enddi
where $\gamma(|)=0$ (the trivial tree has no interior vertices) 
\cite{Butcher_2008}. For instance, the tree factorial of the $n$-th order 
left-comb is 
\begdi
\gamma(\tau_l^{n})=\gamma(|\vee \tau_l^{n-1})=n\gamma(\tau_l^{n-1}).
\enddi
Repeating the procedure $n$ times one arrives at $\gamma(\tau_l^{n})=n!$. Thus, 
the standard factorial is a special case of the tree 
factorial. An analogous procedure applies for right-combs.

The next three lemmas and theorem were developed in order to derive the 
main results of the paper in Section~\ref{sec:4}.  The first lemma provides 
bounds for particular types of non-commutative iterated integrals.

\begle Let $\tau$ be an arbitrary tree in $\mathfrak{T}_n$, 
$\tau_l^{n}$ the left-comb tree in $\mathfrak{T}_n$, $x_i\in X$ and $\eta\in 
X^n$ (all words in $X^\ast$ of length $n$). The non-commutative iterated 
integrals satisfy:
\begin{enumerate}
 \item[$i$.] $\norm{E_{{x_i^n}_{\tau}}[u](t)}_1 
\le\frac{\displaystyle \bar{U}_i^{\abs{\tau}}(t)}{\displaystyle \gamma(\tau)}$, 
\\

\item[$ii$.] $\norm{E_{\eta_{\tau_l^{n}}}[u](t)}_1 \le \prod_{j=1}^n 
\frac{\displaystyle \bar{U}_j^{n_j}(t)}{\displaystyle n_j !}$, \\

\end{enumerate}
where $\bar{U}_j(t) := \int_0^t\bar{u}_i(s) \,ds$, 
$\bar{u}_i(s):=\norm{u_j(s)}_1$, $n_j=\abs{\eta}_{x_j}$ for 
$j=0,\hdots,m$, and $\abs{\eta}_{x_j}$ denotes the number of $x_j$ letters 
in $\eta\in X^\ast$.
\endle

\begpr Bound $i$ is proved by induction over $n$. The $n=0, 1$ cases are 
trivial. Let $\tau = \tau^1 \vee \tau^2$ with $\tau^1\in \mathfrak{T}_k$ and 
$\tau^2 \in \mathfrak{T}_{n-k-1}$ for $0\le k \le n-1$. Assume $i$ holds for 
any $k < n$. Then
\begin{align*}
\lefteqn{\norm{E_{{x_i^n}_\tau}[u](t)}_1}\\
\le & \int_0^t \norm{E_{{x_i^{\abs{\tau^1}}}_{\tau^1}}[u](s)}_1 
\norm{u_i(s)}_1 \norm{E_{{x_i^{\abs{\tau^2}}}_{\tau^2}}[u](s)}_1 \,ds   \\
\le & \int_0^t  \bar{u}_i(s) \frac{\bar{U}_i^{\abs{\tau^1}}}{\gamma(\tau^1)} 
\frac{\bar{U}_i^{\abs{\tau^2}}}{\gamma(\tau^2)} \,ds   \\
= & \frac{\bar{U}_i^{\abs{\tau^1}+\abs{\tau^2}+1}}{(\abs{\tau^1}+\abs{\tau^2}
+1)\gamma(\tau^1)\gamma(\tau^2) } =  
\frac{\bar{U}_i^{\abs{\tau}}}{\gamma(\tau)}. 
\end{align*}
Thus, $i$ holds for all $n\ge 0$. 

Bound $ii$ is also proved by induction over $\abs{\tau_l^{n}}=n$. The $n=0,1$ 
cases are trivial. Let $\eta=x_i\eta'$ with $\eta'\in X^{n-1}$, and recall 
$\tau_l^{n}= | \vee \tau_l^{n-1}$ with $\tau_l^{n-1}$ the $(n-1)$-th 
left-comb. If $ii$ holds for $n-1$, then 
\begin{align*}
\norm{E_{\eta_{\tau_l^{n}}}[u](t)}_1 \le & \int_0^t \norm{u_i(s)}_1 
\norm{E_{\eta'_{\tau_l^{n-1}}}[u](s)}_1 \,ds  \\
\le & \int_0^t  \bar{u}_i(s) \frac{\bar{U}_1^{n'_1}(s) \cdots 
\bar{U}_m^{n'_m}(s)}{n'_1!\cdots n'_m!}  \,ds   \\
\le & \prod_{\substack{j=1 \\ j\neq i}}^m 
\frac{\bar{U}_j^{n'_j}(t)}{n'_j!}\int_0^t  \bar{u}_i(s) 
\frac{\bar{U}_i^{n'_i}}{n'_i!} \,ds   \\
= & \frac{\bar{U}_i^{n'_i+1}}{(n'_i+1)!}\prod_{\substack{j=1 \\ j\neq i}}^m 
\frac{\bar{U}_j^{n'_j}(t)}{n'_j!}=\prod_{j=1}^m 
\frac{\bar{U}_j^{n_j}(t)}{n_j!}, 
\end{align*}
where $n_j=n'_j+1$ and $n_i=n'_i$ for $i\neq j$. So $ii$ applies for all $n\ge 
0$.
\endpr

The following lemma provides a relationship between commutative and 
non-commutative iterated integrals. It plays a key role in the 
convergence analysis of dendriform Fliess operators. 

\begle \label{le:replace_u_with_baru} Let $\eta_\tau\in \mathfrak{T}X^\ast$, 
and $u\in B^{m\times(n\times n)}_1(R)[0,T]$. The 
iterated integral corresponding to $\eta_\tau$ satisfies
\begin{align*}
\norm{E_{\eta_\tau}[u](t)}_1 \le E_{\eta_\tau}[\bar{u}](t), \,\forall t\in 
[0,T], 
\end{align*}
where $\bar{u}(t)=(\bar{u}_1(t), \hdots, \bar{u}_m(t))^T$.
\endle

\begpr The lemma is proved by induction over $\abs{\eta_\tau}=n$. The result is 
trivial for $n=0$. For $n=1$, 
\begin{align*}
\norm{E_{{x_i}}[u](t)}_1 = & {\,} \max_j \sum_{l=1}^n \abs{\int_{0}^t 
(u_i)_{lj}(s)\,ds} \\
\le & {\,}    \int_{0}^t \max_j \sum_{l=1}^n\abs{(u_i)_{lj}(s)}\,ds \\
= & {\,}    \int_{0}^t \norm{u_i(s)}_1\,ds \\
= & {\,}    \int_{0}^t \bar{u}_i(s) \,ds = E_{x_i}[\bar{u}](t),
\end{align*}
where $\bar{u}_{i}(t) \ge 0 $ is now scalar-valued, i.e., is a commutative 
input. If the claim now holds up to order $n$ and $\eta_\tau= 
\Phi^{-1}(\tau^1_\xi 
\vee_{x_i} \tau^2_{\nu})$, then
\begin{align*}
\lefteqn{\hspace*{-0.3in}\norm{E_{\eta_\tau}[u](t)}_1} \\
\hspace*{-0.2in}\le & {\,} \int_{0}^t \norm{E_{\xi_{\tau^1}}[u](s)}_1 
\norm{u_i(s)}_1 
\norm{E_{\nu_{\tau^2}}[u](s)}_1 \,ds \\
\le & {\,} \int_{0}^t E_{\xi_{\tau^1}}[\bar{u}](s) \bar{u}_i(s) 
E_{\nu_{\tau^2}}[\bar{u}](s) \,ds \\
= & {\,} E_{\eta_\tau}[\bar{u}](t).
\end{align*}
Thus, the bound holds for all $n\ge 0$.
\endpr

It is important to note that even though the components of $\bar{u}$ are 
mutually commutative, the corresponding iterated integrals do not coincide with 
the commutative counterpart where one removes the ordering provided by 
the trees. 
\begex Let $\eta=x_ix_jx_k$ ($i\neq j \neq k$) and 
$\tau=\arbsmall{131}$. 
Then it follows that 
\begin{align*}
\lefteqn{\hspace*{-0.4in} E_{(x_ix_jx_k)_\tau}[\bar{u}](t)} \\
\hspace*{0.4in} = & {\,} \int_0^t \left(\int_0^s \bar{u}_i(r)\, dr\right) 
\bar{u}_j(s) \left(\int_0^s \bar{u}_k(r)\, dr\right)\,ds \\
= & {\,} \int_0^t \bar{u}_j(s)\left(\int_0^s \bar{u}_i(r)\, dr\right)  
\left(\int_0^s \bar{u}_k(r)\, dr\right)\,ds \\
= & {\,} E_{x_jx_ix_k}[\bar{u}](t)+E_{x_jx_kx_i}[\bar{u}](t),
\end{align*}
where $E_{x_jx_ix_k}[\bar{u}](t)$ and $E_{x_jx_kx_i}[\bar{u}](t)$ are 
commutative iterated integrals distinct from $E_\eta[\bar{u}](t)$.
\endex

The correspondence between the commutative shuffle product, $\shuffle$, and the 
product of commutative iterated integrals generalizes in the non-commutative 
setting as follows.  
\begth \label{th:shuffleNC_product_theorem}
Let $u\in B^{m\times(n\times n)}_1(R)[0,T]$ and 
${\eta}_{\tau^1},{\xi}_{\tau^2}\in \mathfrak{T}X^\ast$. Then
\begeq \label{eq:non-commut_iterated_integrals_shuffle}
E_{\eta_{\tau^1}}[u](t) E_{\xi_{\tau^2}}[u](t) = 
E_{\eta_{\tau^1}\shuffleNC \xi_{\tau^2}}[u](t).
\endeq
\endth

\begpr
Recall that the decorated tree corresponding to $\eta_{\tau^i}$ is 
$\tau^i_{\eta}=\Phi(\eta_{\tau^i})$. Identity 
\rref{eq:non-commut_iterated_integrals_shuffle} is proved by 
induction over $|\tau^1_\eta|+|\tau^2_\xi|=n$. The claim is 
trivial for $n=0,1$ since $E_{\emptyset}[u]=I$ and by definition $\eta_{\tau^1} 
\shuffleNC \emptyset =  \emptyset \shuffleNC \eta_{\tau^1}=\eta_{\tau^1}$. 
Assume \rref{eq:non-commut_iterated_integrals_shuffle} holds up to some 
fixed $n\ge 1$. If $\tau^{1}_\eta=\tau^{11}_{\eta_1}\vee_{x_i} 
\tau^{12}_{\eta_2}$ and 
$\tau^{2}_{\xi}=\tau^{21}_{\xi_1}\vee_{x_j}\tau^{22}_{\xi_2}$ with 
$|\tau^1_\eta|+|\tau^2_{\xi}|=n+1$, then 
\begin{align*}
\lefteqn{\hspace*{-0.2in} E_{\eta_{\tau^1}}[u](t)  E_{\xi_{\tau^2}}[u](t)}\\
= & {\,} \int_0^t E_{\eta_{\tau^{11}}}[u](s) u_i(s) E_{\eta_{\tau^{12}}}[u](s) 
\, ds  \\
&  \int_0^t E_{\xi_{\tau^{21}}}[u](s)u_j(s)E_{\xi_{\tau^{22}}}[u](s) \,ds \\
= & {\,} \int_0^t E_{\eta_{\tau^{11}}}[u](s) u_i(s) E_{\eta_{\tau^{12}}}[u](s) 
\\
  & \left(\int_0^s E_{\xi_{\tau^{21}}}[u](r) u_j(r) E_{\xi_{\tau^{22}}}[u](r) 
\,dr \right) ds \\
& + \int_0^t \left(\int_0^s E_{\eta_{\tau^{11}}}[u](r) u_i(r) 
E_{\eta_{\tau^{12}}}[u](r) \,dr \right) \\
& E_{\xi_{\tau^{21}}}[u](s) u_j(s) E_{\xi_{\tau^{22}}}[u](s)\,ds \\
= & {\,} \int_0^t E_{\eta_{\tau^{11}}}[u](s) u_i(s) E_{\eta_{\tau^{12}} 
\shuffleNC \xi_{\tau^2}}[u](s)\, ds  \\
& + \int_0^t E_{\eta_{\tau^1} \shuffleNC 
\xi_{\tau^{21}}}[u](s) u_j(s)E_{\xi_{\tau^{22}}}[u](s)\,ds \\
= & {\,} E_{\Phi^{-1}(\underbrace{\text{\scriptsize 
$\tau^{11}_{\eta_1}\vee_{x_i}(\tau^{12}_{\eta_2} \shuffle 
\tau^2_{\xi})$}}_{\text{\scriptsize $\tau^1_{\eta}\prec 
\tau^2_{\xi}$}})}[u](t)\\
  & + E_{\Phi^{-1}(\underbrace{\text{\scriptsize $(\tau^{1}_{\eta} \shuffle 
\tau^{21}_{\xi_1})\vee_{x_j} \tau^{22}_{\xi_2}$}}_{\text{\scriptsize 
$\tau^1_{\eta} \succ \tau^2_{\xi}$}})}[u](t)  \\
= & {\,} E_{\eta_{\tau^1} \prec \xi_{\tau^2}}[u](t) + E_{\eta_{\tau^1} \succ 
\xi_{\tau^2}}[u](t) =  E_{\eta_{\tau^1} \shuffleNC \xi_{\tau^2}}[u](t).
\end{align*}
\endpr

The final lemma in the section is the result of the 
grouping of trees with same order (Lemma 
\ref{le:characteristic_shuffle_identity}) and Lemma 
\ref{le:replace_u_with_baru}. 

\begle \label{le:non-commutative_shuffle identity_for ubar}
Let $\tau=\arbsmall{1}$. The following identity holds when $u_i$ is replaced 
with $\bar{u}_i$:
\begdi %eq \label{eq:non-commutative_shuffle identity_for ubar}
E_{{(x_i)}^{\shuffleNC n}_\tau}[\bar{u}](t)= n! E_{x_i^n}[\bar{u}](t).
\enddi
\endle
\vspace*{0.05in}

\begpr For brevity define $x_i^{\shuffleNC n} = {(x_i)}^{\shuffleNC n}_\tau$ 
and 
recall that $ x_i^{\shuffleNC n} = \Phi^{-1}\left(\arbsmall{1}^{\shuffleNC 
n}\right) \in \allpolyTX$. In the commutative setting, the definition of an 
iterated integral coincides with the ordering of a non-commutative iterated 
integral corresponding to left-comb trees (see \rref{eq:comm_E_map}). Thus, 
replacing $u_i$ with $\bar{u}_i$, one has
\begin{align*}
E_{x_i^{\shuffleNC n}}[u](t) = E_{x_i^{\shuffle n}}[\bar{u}](t).
\end{align*}
Applying the identity $x_i^{\shuffle n}=n! x_i^n$ proves the lemma.
\endpr

\section{Dendriform Fliess operators and their convergence} \label{sec:4}

In this section dendriform Fliess operators are defined, and 
sufficient conditions for their convergence are provided.

\subsection{Dendriform Fliess operators}

The definition of a dendriform Fliess operator is given first.
\begde \label{de:Dendriform_Fliess_operator}
Let $u\in B^{m\times (n \times n)}_1(R)[0,T]$ and $c\in \allseriesTXelln$. 
A \emph{dendriform Fliess operator} with generating series $c$ is defined 
by the following summation
\begeq \label{eq:Dendriform_Fliess_operator}
F_c[u](t) = \sum_{\eta_\tau \in \mathfrak{T}X^\ast} (c,\eta_\tau) 
E_{\eta_\tau}[u](t).
\endeq
%which is referred as a \emph{dendriform Fliess operator}.
\endde

The operator in \rref{eq:dyson_series} is a special case of a dendriform 
Fliess operator. The support of its generating series contains only left-comb 
trees. This is purely a consequence of the iterative procedure used to derive 
it. However, defining Fliess operators as a summation comprised of only 
left-comb trees limits its application as shown in the next example.
\vspace*{-0.2in}
\begin{figure}[H]
\vspace*{-0.1in}
\caption{Product connection of Fliess operators}
\label{fig:product_connection}
\begce
\includegraphics[width=5cm]{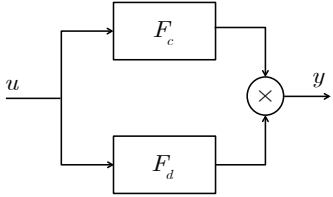}
\endce
\vspace*{-0.1in}
\end{figure}
\vspace*{-0.1in}
\begex \label{ex:product_connection}
Suppose two Fliess operators $F_c$ and $F_d$ have generating series 
in terms of left-combs. Assume $c=c'I$ 
and $d=d'I$ in $\allseriesTXnn$ with $c',d'$ being scalar-valued series, and 
$u\in B^{1\times (n \times n)}_1(R)[0,T]$ for some $R,T > 0$. Since 
$\sum_{\eta_\tau\in \mathfrak{T}X^\ast} \eta_\tau = \sum_{\eta\in X^\ast} 
\sum_{\tau\in \mathfrak{T}} \eta_\tau $, their product 
connection as shown in Figure \ref{fig:product_connection} is described by
\begin{align*}
F_c[u] F_d[u]& = \!\sum_{\hspace*{-0.1in}n_1,n_2=0}^\infty 
\!\!\!\!\sum_{\;\;\;\substack{\eta\in X^{^{n_2}}\\ \xi 
\in X^{^{n_1}}}} \!\!(c,\eta_{\tau_l^{n_1}}) 
(d,\xi_{\tau_l^{n_2}}) E_{\eta_{\tau_l^{n_1}}}[u] E_{\xi_{\tau_l^{n_2}}}[u].
\end{align*}
Recall that $E_{\eta_{\tau_l^{n_1}}}[u] E_{\xi_{\tau_l^{n_2}}}[u] = 
E_{\eta_{\tau_l^{n_1}}\shuffleNC \xi_{\tau_l^{n_2}} }[u]$, where $\shuffleNC$ 
generates \emph{more} than just left-combs as shown in 
Example~\ref{ex:non-commutative_shuffle}. Therefore, Definition 
\ref{de:Dendriform_Fliess_operator} is general enough to characterize  
such interconnections in the non-commutative framework.
\endex

\subsection{Convergence of dendriform Fliess operators}

The next theorem addresses the convergence of dendriform Fliess 
operators by considering bounds on the coefficients of the corresponding 
generating series. The final three lemmas and 
theorem in Section \ref{sec:3} were specifically developed for proving this 
theorem.

\begth  \label{th:convergence_Fliess_operator}
Let $c\in \allseriesTXelln$ with coefficients satisfying the growth condition 
\begeq \label{eq:coefficients_Dendriform_Fliess_op}
\norm{(c,\eta_\tau)}_1\le K M^{\abs{\tau}}, \quad \forall \eta_\tau \in 
\mathfrak{T}X^\ast
\endeq
for some constants $K,M>0$. Then there exist $R,T>0$ such that for each $u \in 
B^{m\times(n\times 
n)}_1(R)[0,T]$ the series
\begdi
y(t)=F_{c}[u](t)= \sum_{\eta_\tau \in \mathfrak{T}X^\ast} (c,\eta_\tau) 
E_{\eta_\tau}[u](t)
\enddi
converges absolutely and uniformly on $[0,T]$.
\endth

\begpr
Fix some $T>0$. Pick $u\in B^{m\times(n\times n)}_1(R)[0,T]$ and let $R:= 
\max\{\norm{u},T\}$. Since the summation over dendriform words can be 
decomposed into the summations over words in $X^\ast$ (decorations) and the 
summation over trees, define
\begdi
a_k = \sum_{\eta\in X^k} \sum_{\tau\in \mathfrak{T}_k} (c,\eta_\tau) 
E_{\eta_\tau}[u].
\enddi
Using \rref{eq:coefficients_Dendriform_Fliess_op} and Lemma 
\ref{le:replace_u_with_baru}, a bound for $a_k(t)$ is computed as
\begin{align*} 
\nonumber \norm{a_k} = & {\,} \norm{ \sum_{\eta\in X^k} 
\sum_{\tau\in \mathfrak{T}_k} (c,\eta_\tau) E_{\eta_\tau}[u]} \\
\nonumber \le & {\,}   \sum_{\eta\in X^k} \norm{(c,\eta_\tau)} \sum_{\tau\in 
\mathfrak{T}_k}  \norm{E_{\eta_\tau}[u]} \\
\le & {\,}   KM^k \sum_{\eta\in X^k}  \sum_{\tau\in \mathfrak{T}_k}  
E_{\eta_\tau}[\bar{u}].
\end{align*}
From \rref{eq:non-commut_iterated_integrals_shuffle}, Lemma 
\ref{le:characteristic_shuffle_identity} and the commutativity of 
$\bar{u}$, one has that
\begin{align*}
\lefteqn{\hspace*{-0.4in}\sum_{\eta\in X^k} \sum_{\tau\in \mathfrak{T}_k}  
E_{\eta_\tau}[\bar{u}] = E_{\Phi^{-1}\left(\sum_{\eta\in X^k} \eta; 
\sum_{\tau\in\mathfrak{T}_k} \tau\right)} [\bar{u}]} \\
\hspace*{0.4in} = & \sum_{\alpha_0+\cdots 
+\alpha_m=k}E_{\Phi^{-1}\left(x_1^{\alpha_0}\shuffle 
\cdots \shuffle x_m^{\alpha_m}; \arbsmall{1}^{\shuffleNC k} 
\right)}[\bar{u}].
\end{align*}
Lemma \ref{le:non-commutative_shuffle identity_for ubar} in tree terminology 
amounts to $(\tau_l^1)^{\shuffleNC k}= k! \tau_l^k $. This is also equivalent to
\begin{align*}
E_{\Phi^{^{-1}}\left(x_1^{\alpha_0}\shuffle \cdots  
\shuffle x_m^{\alpha_m}; \arbsmall{1}^{\shuffleNC k} \right)}[\bar{u}]
= & {\,} k! E_{ x_1^{\alpha_0}\shuffle \cdots \shuffle 
x_m^{\alpha_m}}[\bar{u}].
\end{align*}
Continuing the analysis,
\begin{align*}
\lefteqn{\hspace*{-0.7in}\sum_{\eta\in X^k} \sum_{\tau\in \mathfrak{T}_k}  
E_{\eta_\tau}[\bar{u}]=  k! \!\!\!\!\! \sum_{\alpha_0+\cdots +\alpha_m = k}E_{ 
x_1^{\alpha_0}\shuffle \cdots \shuffle x_m^{\alpha_m}}[\bar{u}]} \\
\hspace*{0.4in} = & {\,} \sum_{\alpha_0+\cdots +\alpha_m = k} 
\frac{k!}{\alpha_0!\cdots 
\alpha_m! } E_{x_1}^{\alpha_0}[\bar{u}] \cdots E_{x_m}^{\alpha_m}[\bar{u}]\\
\le &  {\,\;} R^k\!\!\!\!\!\sum_{\alpha_0+\cdots +\alpha_m = k} 
\frac{k!}{\alpha_0!\cdots 
\alpha_m! } = ((m+1)R)^k,
\end{align*}
where $E_{x_i}[\bar{u}](t) = \bar{U}_i(t) \le \norm{u} \le R$. It is now clear 
that 
\begin{align*}
\sum_{k= 0}^\infty \norm{a_k(t)} \le {\,} &  \sum_{k= 0}^\infty K (M R (m+1))^k.
\end{align*}
Therefore, $F_c[u](t)$ converges absolutely and uniformly on $[0,T]$ for $R < 
\frac{1}{M(m+1)}$.
\endpr

Coefficients bounded as in \rref{eq:coefficients_Dendriform_Fliess_op} give 
convergence of a local nature whereas in the commutative case such 
coefficients bounds provide a type of global convergence \cite{Gray-Wang_SCL02}. 
The reason for this discrepancy is that in addition to summing over all possible 
permutations of letters in $X$, which is the commutative case, the bounds for 
non-commutative iterated integrals also require the summation over all trees. 
This contributes an extra $k!$ factor coming directly from the integrals.

A \emph{left-comb dendriform Fliess operator} is a dendriform Fliess 
operator whose generating series support only have dendriform 
words corresponding to left-combs. The convergence of such operators is 
addressed in the next theorem. 
\begth \label{th:convergence_Fliess_operator_left_comb}
Let $c\in \allseriesTXelln$ with coefficients satisfying the growth condition 
\begdi
\norm{(c,\eta_{\tau})}_1 \le KM^{\abs{\tau}}\abs{\tau}!, \quad \forall 
\eta_\tau\in \mathfrak{T}X^\ast
\enddi
for some constants $K,M>0$ and $\supp(c)\subseteq \{\eta_\tau \in 
\mathfrak{T}X^\ast, \tau = \tau^k_l,k>0\}$. Then 
there exist $R,T>0$ such that for each $u \in B^{m\times(n\times 
n)}_1(R)[0,T]$ the series
\begeq \label{eq:left-comb_Fliess_operator}
y(t)=F_{c}[u](t)= \sum_{k= 0}^\infty \sum_{\eta \in X^k} (c,\eta_{\tau^k}) 
E_{\eta_{\tau^k}}[u](t)
\endeq
converges absolutely and uniformly on $[0,T]$.
\endth

\begpr
The proof is similar to the one for Theorem 
\ref{th:convergence_Fliess_operator}. However, there is no $k!$ factor from the 
iterated integrals since the series only depends on left-combs.
\endpr

\begex Consider $c = \sum_{k= 0}^\infty {x_1^k}_{\tau_l^k}\in \allseriesTXnn$. 
This series is the generating series corresponding to \rref{eq:dyson_series}, 
which is the solution of \rref{eq:operator_evolution}. Recall that 
\rref{eq:operator_evolution} can represent the evolution of a closed quantum 
system (all quantum constants normalized to $1$). In the commutative case, it is 
known that $X =\exp(\Omega)$, where $\Omega(t)=\int_0^t U(s)\,ds$. From the 
Fliess operator point of view, 
\begdi 
Z = F_c[U] = \sum_{n= 0}^\infty  E_{x_1^n}[U],
\enddi
which by the properties of the commutative shuffle product gives
\begeq \label{eq:example_NC_x1star}
Z=F_c[U] = \sum_{k= 0}^\infty  \frac{\left(E_{x_1}[U]\right)^k}{k!} = 
\exp\left({E_{x_1}[U]}\right),
\endeq
where obviously $ E_{x_1}[U]=\Omega$. Suppose now $U$ is non-commutative. 
Then
\begin{align} \label{eq:example_NC_x1star2}
F_c[U] & = \sum_{k= 0}^\infty  E_{{x_1^k}_{\tau_l^k}}[U],
\end{align}
which by Theorem \ref{th:convergence_Fliess_operator_left_comb} with $K=M=1$ 
is well defined. Assume now that $F_c[U]$ has an exponential representation 
similar to the commutative case. That is, $F_c[U] = \exp(\Omega)$ with $\Omega= 
F_d[U]$ for some $d\in\allseriesTXnn$. Unfortunately, the identities used to 
obtain \rref{eq:example_NC_x1star} cannot be used to find the expression for 
$d$. But Lemma \ref{le:characteristic_shuffle_identity} 
provides an inductive way to compute it. Assume that $\Omega = 
E_{{x_1}_{\tau_l^1}}[U]$. Then expanding $\exp(\Omega)$ gives
\begin{align*}
\exp\left({\Omega}\right)= & {\,} I + E_{{x_1}_{\tau_l^1}}[U] + 
\frac{1}{2!}\left(E_{{x_1}_{\tau_l^1}}[U]\right)^2 + \cdots\\
 = & {\,} I + E_{{x_1}_{\tau_l^1}}[U] + 
\frac{1}{2!} E_{\Phi^{^{-1}}\left(x_1^2;\arbsmall{12} + 
\arbsmall{21}\right)}[U]+\cdots\\
 = & {\,} I + E_{{x_1}_{\tau_l^1}}[U]+
E_{{x_1^2}_{\tau_l^2}}[U] \\
 & -\!\frac{1}{2!} E_{\Phi^{^{-1}}\!\left(x_1^2;\arbsmall{21}\right)}[U] + 
\frac{1}{2!} 
E_{\Phi^{^{-1}}\!\left(x_1^2;\arbsmall{12}\right)}[U] +\!\cdots
\end{align*}
Observe that the expansion produces more terms than needed. Therefore, a
correction term must be used in order cancel the extra second order terms. So 
redefine $\Omega$ as 
\begdi 
\Omega = E_{{x_1}_{\tau_l^1}}[U]-\frac{1}{2} E_{\Phi^{^{-1}}\left( x_1^2; 
\arbsmall{12} - \arbsmall{21}\right)}[u].
\enddi 
It follows then that the first and second order terms are
\begin{align*}
\exp\left({\Omega}\right) = & {\,} I + E_{{x_1}_{\tau_l^1}}[U] - 
\frac{1}{2}E_{\Phi^{^{-1}}\left(x_1^2;\arbsmall{12} - 
\arbsmall{21}\right)}[U] + \\
& +\! \frac{1}{2!} E_{\Phi^{^{-1}}\!\left(x_1^2;\arbsmall{12}\right)}[U] + 
\frac{1}{2!} E_{\Phi^{^{-1}}\!\left(x_1^2;\arbsmall{21}\right)}[U] +\! \cdots \\
 = & {\,}  I + E_{{x_1}_{\tau_l^1}}[U] + E_{{x_1^2}_{\tau_l^2}}[U] + \cdots
\end{align*}
Observe that the second order correction term in $\Omega$ can be written in 
the following form
\begin{align*}
E_{\Phi^{^{-1}}\left( x_1^2; 
\arbsmall{12}-\arbsmall{21}\right)}[U] 
= E_{x_1\prec x_1}[U] - E_{x_1\succ x_1}[U].
\end{align*}
In fact, defining the product $\triangleright = \prec - \succ$, it 
follows that 
\begin{align*}
E_{x_1 \triangleright x_1}[U](t) & = E_{x_1\prec x_1}[U](t) - 
E_{x_1\succ x_1}[U](t) \\
& =   \int_{0}^t \left[U(s),\int_0^s U(r)\,dr\right]ds,
\end{align*}
where $[\cdot,\cdot]$ representing the commutator. The non-associative product 
$\triangleright$ is an example of a \emph{pre-Lie 
product} \cite{Kurusch-Manchon_2009}. This correction procedure can 
be applied successively at every order. At order $3$, the correction terms for 
$\Omega$ are
\begin{align*}
\Omega = & {\,}  E_{x_1}[U] -\frac{1}{2} E_{x_1 \triangleright x_1}[U] 
+ \frac{1}{4} E_{(x_1 \triangleright x_1)\triangleright x_1}[U] \\
& + \frac{1}{12} E_{x_1 \triangleright (x_1 \triangleright x_1)}[U],	
\end{align*}
which gives
\begin{align*}
\exp\left({\Omega}\right) \!= \!  I + 
E_{{x_1}_{\tau_l^1}}[U] + E_{ {x_1^2}_{\tau_l^2}}[U] + 
E_{{x_1^3}_{\tau_l^3}}[U] 
+\cdots\!,
\end{align*}
where the generating series $d$ of $\Omega$ satisfies the recursion 
\begin{align*}
d^{[k]} = & {\,} \sum_{n\ge 0} \frac{B_n}{n!} L^{(n)}_{d^{[k-1]} 
\triangleright} (x_1)
\end{align*}
with $d^{[1]}=x_1$, $\lim_{k \rightarrow \infty} d^{[k]}= d$, $L^{(n)}_{d 
\triangleright}(x) = d \triangleright (L^{(n-1)}_{d \triangleright}(x))$,  
$L^{(0)}_{d \triangleright}(x)=d$, and $B_n$ denotes the $n$-th Bernoulli 
number. Thus, the limit of $\exp(F_d^{[k]}[U])$ as $k\rightarrow 
\infty$ agrees with \rref{eq:example_NC_x1star2}. This is the 
well-known \emph{Magnus expansion}. The more familiar expression for the 
Magnus expansion is obtained by noting that
\begin{align*}
E_{L^{(n)}_{d \triangleright}(x_1)}[U](t)  & = \int_0^t 
ad^{(n)}_{\Omega(s)}(U(s))\,ds, 
\end{align*}
and $ad^{(n)}_{\Omega}(U)=[\Omega,ad^{(n-1)}_{\Omega}(U)]$ with 
$ad^{(0)}_{\Omega}(U)=U$. Compared to the ordered exponential presented in the 
introduction, this is a true exponential. In quantum mechanics this is one way 
to show that the evolution operator is unitary for all times. Finally, the 
Fliess operator $F_c[U]$ in \rref{eq:example_NC_x1star2} provides an 
input-output map that encodes in the iterated integrals the underlying algebraic 
structure of the system. 
\endex

\section{Conclusions} \label{sec:conclusions}

A setting for dendriform Fliess operators has been provided. The algebraic 
structure basically considers the relationship between dendriform words and 
trees. Sufficient conditions for the convergence of such Fliess operators were 
given for the general case \rref{eq:Dendriform_Fliess_operator} and for 
operators indexed only by left-comb trees \rref{eq:left-comb_Fliess_operator}.

\end{document}